\documentclass[10pt,a4paper]{article}

\usepackage{floatflt,graphicx,graphics}
\usepackage{amsmath,amssymb,color}
\usepackage{amsthm}
\usepackage[autolinebreaks]{mcode}
\usepackage{url}

\newcommand{\C}{\mathbb{C}}

\newcommand{\D}{\mathbb{D}}
\renewcommand{\i}{\mathrm{i}}
\newcommand{\bI}{{\bf I}}
\newcommand{\bM}{{\bf M}}
\newcommand{\bN}{{\bf N}}
\newcommand{\capa}{\mathrm{cap}}
\newcommand{\caph}{\mathrm{caph}}
\newcommand{\cape}{\mathrm{cape}}
\newcommand{\cc}{\mathrm{c}}
\renewcommand{\Im}{{ \rm Im}\,}
\renewcommand{\Re}{{ \rm Re}\,}

\newtheorem{thm}{Theorem}

\begin{document}

\title{Computation of conformal invariants}

\author{Mohamed M. S. Nasser$^{\rm a}$ and Matti Vuorinen$^{\rm b}$}

\date{}
\maketitle

\vskip-0.5cm %
\centerline{$^{\rm a}$Department of Mathematics, Statistics and Physics, Qatar University, Doha, Qatar.}%
\centerline{mms.nasser@qu.edu.qa}

\vskip0.5cm %
\centerline{$^{\rm b}$Department of Mathematics and Statistics, University of Turku, Turku, Finland.} %
\centerline{vuorinen@utu.fi}

\begin{center}
\begin{quotation}
{\noindent {\bf Abstract.\;\;}%
We study numerical computation of conformal invariants of domains in the complex plane. In particular, we provide an algorithm for computing the conformal capacity of a condenser. The algorithm applies to a wide variety of geometries: domains are assumed to have smooth or piecewise smooth boundaries. The method we use is based on the boundary integral equation method developed and implemented in~\cite{Nas-ETNA}. A characteristic feature of this method is that, with small changes in the code, a wide spectrum of problems can be treated.
We compare the performance and accuracy to previous results in the cases when numerical data is available and also in the case of several model problems where exact results are available.
}%
\end{quotation}
\end{center}
\begin{center}
\begin{quotation}
\end{quotation}
\end{center}

\begin{center}
\begin{quotation}
\end{quotation}
\end{center}
%

\section{Introduction} \label{section1}

During the past fifty years conformal invariants have become crucial tools for complex analysis. Most important of these invariants
are the conformal capacity, the harmonic measure, the extremal length, and the hyperbolic distance \cite{ah,garmar,kela,st}. But this is not all: the generalized capacity, the transfinite diameter, the reduced extremal length, the hyperbolic area, and the modulus metric \cite{du,kir,Vas02,vu85,vu88} are some additional examples, see~\cite[Ch 10]{hkv}. Some of the many applications of these tools are discussed in the articles of the handbook \cite{kuh3}. In view of the plenitude of these applications, it is surprising that these invariants can be expressed explicitly only in very few special cases. Sometimes rudimentary upper or lower bounds for conformal invariants in terms of less involved comparison functions can provide important steps in proofs.

At the same time it seems that the full power of conformal invariance remains unused.
One reason for this is that the analytic expressions for conformal invariants are usually
too complicated for pen and paper calculations and the existing computational methods are
scattered throughout the mathematical literature: the way from theory to practical
experimentation is too long. On the other hand, the creators of the existing computational methods may not be aware of the scope of applicability of their methods in theoretical
studies: the way from experiments to theory is also long. If the distance from theory to
experiments could be made shorter, a theoretical researcher could easily experiment with the dependence of a problem under perturbation of geometry and vice versa a computational scientist would find new types of benchmark problems and areas of application.

With the above ideas as our guiding principles, we have written a series of papers of which the present one is devoted to doubly connected domains~\cite{Nvm,Nvs}.
As far as we know, our work is the first attempt to provide computational tools for a wide class of conformal invariants with the feature that modification of geometry is simple. The method we use was
 developed and implemented by the first author \cite{Nas-ETNA} and we apply it to study
several computational problems never studied before and we also compare its performance
to several results in the literature.
As test problems we use the computation of condenser capacity, a topic studied by the second author in several papers \cite{hrv1,hrv2,hrv3}.

A condenser is a pair $(D,C)$ where $D$ is an open set in $\mathbb{R}^2$ and $C \subset D$ is compact.
In our study we assume that the topology is simple but still general enough for most applications:
the sets $\partial D$ and $\partial C$ are connected sets, each set is a piecewise smooth Jordan curve. The cases for which the sets $\partial D$ and $\partial C$ are slits will be also considered.  

The conformal capacity, or capacity for short, is defined by \cite{du}
\begin{equation} \label{capdef}
{\rm{cap}}(D,C) = \int_D |\nabla u|^2 \, dm
\end{equation}
where $u:D \to \mathbb{R}$ is a harmonic function with $u(x) \ge 1\,$ for all $x \in C\,$  and $u(x) \to 0\,$  for  $x \to \partial D\,.$
The domain $G=D\setminus C$ is called the field of the condenser and the closed sets $C$ and $D^c$ are called the plates of the condenser. Then, the capacity ${\rm{cap}}(D,C)$ may alternatively be written as ${\rm{cap}}(G)$.

In literature, only very few formulas are given for the capacities of concrete condensers.
Numerical methods are therefore needed to compute the value of~\eqref{capdef}.
Our problem is reduced to the classical problem of solving numerically the Dirichlet boundary value problem for the Laplace equation. Moreover, by the Dirichlet principle \cite[pp. 447-456]{garmar},
the extremal function $u_0$ is harmonic and minimizes the integral \cite{ah}, \cite[pp.441-456]{garmar}:
\begin{equation} \label{capdef2}
 \int_D |\nabla u_0|^2 \, dm = \inf \{ \int_D |\nabla u|^2 \, dm \}
\end{equation}
where the infimum is taken over all $C_o(D)$ functions with the indicated Dirichlet boundary values.
The capacity of condensers is invariant under conformal mappings, and hence domains with difficult geometry can be treated using conformal mappings~\cite{avv,dt,du,Kyt,ps,sl,Vas02,vu88}. See also~\cite{bbgghv,dek,et}.

Before proceeding to the contents of our work a few general remarks about the literature we know about may be in order. Because of the wide scope of conformal invariants, relatively few cases exist where ``the right answer'' is known. In cases for which the analytic formulas are unknown, the computational performance may be analysed by observing convergence features of the results under successive refinements  of the numerical model, and error estimates maybe based on general
theory. In those relatively few cases we have found in the literature where the analytic
formula is known, the true error estimate may be given. Sometimes a high accuracy can be
achieved, say 12 decimal places, but the dilemma is that if the geometry of the problem is smoothly changed a bit, the method might not be applicable at all.

Section 2 summarizes our computational workhorse, the boundary integral method geared for the capacity computation of ring domains, which will be applied in several later sections, sometimes together with auxiliary procedures. 
In Section 3 we consider ring domains for which the exact value of capacity is known and investigate the performance of our method.
Sections 4, 5, and 6 deal with condensers whose one or two complementary components are slits---these are well-known examples of computationally challenging problems and we use here auxiliary conformal mapping to overcome computational difficulties. 
Section 7 deals with the case when both complementary components of a ring domain are thin rectangles. 
In Section 8, we consider the numerical computation of the hyperbolic capacity and the elliptic capacity of compact and closed sets. 
The final Section 9 gives some concluding remarks and information about the access to our MATLAB software.

\section{Conformal mapping onto annulus}
\label{sec-cm}

\subsection{Ring domains.}\label{ringDef}

A  domain $G$ in the extended complex plane $\overline{\C}=\C\cup\{\infty\}\,,$ whose complement
 $\overline{\C} \setminus G$ has two components, is called a ring domain or, briefly, a ring. It is
 a classical fact that a ring can be mapped by a conformal map onto an annulus $\{   z : q<|z|<1 \}\,, q \in (0,1)\,.$
A ring $R$ is the simplest example of a condenser and its capacity is given by \cite{du}, \cite[p. 132-133]{garmar}
$$
{\rm cap} (R) = \frac{2 \pi}{ \log (1/q)}\,.
$$
The number $ \log (1/q)$ is called the modulus of the ring, i.e.,
\begin{equation}\label{eq:Mod-R}
M(R)=\log(1/q)=\frac{2\pi}{\capa(R)} \,.
\end{equation}
Because of the conformal invariance of the capacity, this definition is independent of the conformal map. For the computation of the capacity we will often use an auxiliary
conformal mapping to avoid computational singularities.

In this section we describe the method of our numerical work, based on the solution of the boundary integral equation with the generalized Neumann kernel~\cite{Nas-ETNA,Weg-Nas}. The integral equation has been applied to calculate conformal mappings onto several canonical domains~\cite{Nas-Siam1,N18,NF13}. We review the application of the integral equation to compute the conformal mapping from doubly connected domains onto an annulus $\{z : q<|z|<1 \}, q \in (0,1)$, and present the MATLAB implementation of the method.
In later sections we will apply this method
for capacity computation of several condensers, in particular, we will consider several types of rings with a simple
geometric structure.

\subsection{The generalized Neumann kernel.}\label{sc:gnk}

\begin{figure}[hbt] %
\centerline{
\hfill
\scalebox{0.4}{\includegraphics[trim=0 0 0 0,clip]{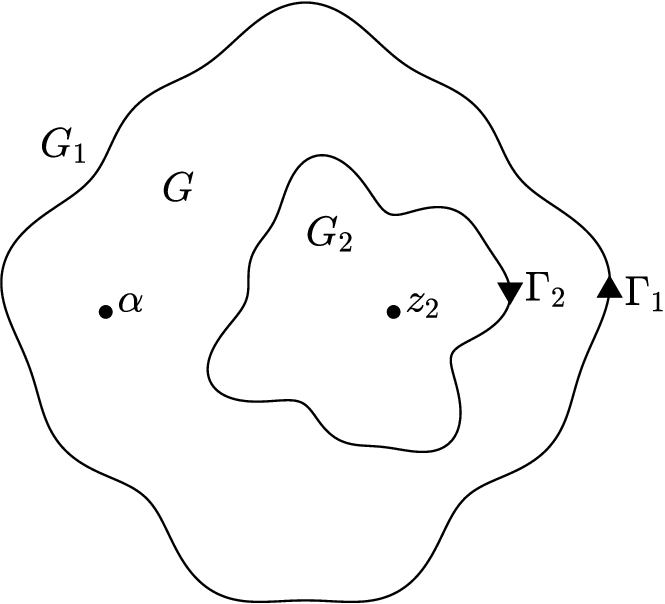}}
\hfill
\scalebox{0.4}{\includegraphics[trim=0 0 0 0,clip]{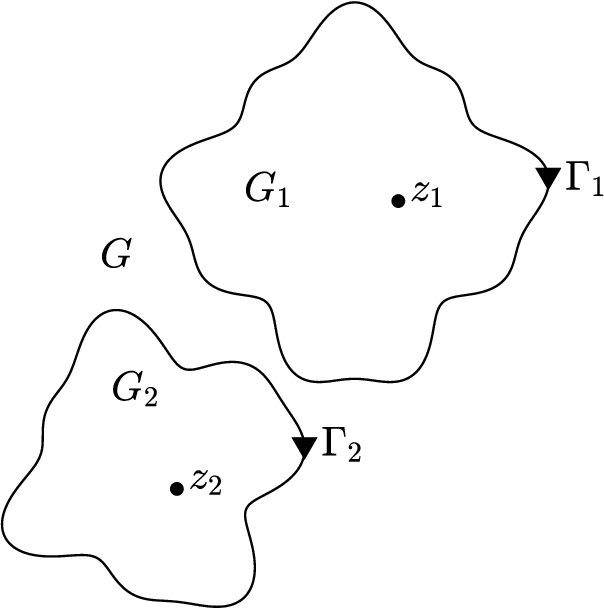}}
\hfill
}
\caption{An example of a bounded doubly connected domain (left) and an unbounded doubly connected domain (right).}
\label{fig:domain}
\end{figure}

Let $G$ be a bounded or an unbounded doubly connected domain bordered by
\[
\Gamma=\partial G=\Gamma_1\cup\Gamma_2
\]
where each of the boundary components $\Gamma_1$ and $\Gamma_2$ is a closed smooth Jordan curve. We choose the orientation of boundary $\Gamma$ such that when we proceed along $ \Gamma\,,$ the domain $G$ is always on the left side. If $G$ is bounded, then $\Gamma_1$ is the external boundary and $\Gamma_2$ is contained in the bounded domain whose boundary is $\Gamma_1$.
The complement $\overline{\C}\backslash\overline{G}$ of the domain $G$ with respect to the extended complex plane $\overline{\C}$ consists of two simply connected domains $G_1$ on the right of $\Gamma_1$ and $G_2$ on the right of $\Gamma_2$. 
The domain $G_2$ is bounded and the domain $G_1$ is unbounded with $\infty\in G_1$. Further, we assume that $\alpha$ is an auxiliary given point in the domain $G$ and $z_2$ is an auxiliary given point in the simply connected domain $G_2$.
When $G$ is unbounded, then $\infty\in G$ and the two domains $G_1$ and $G_2$ are bounded. We assume that $z_1$ and $z_2$ are auxiliary given points in the simply connected domains $G_1$ and $G_2$, respectively.
See Figure~\ref{fig:domain}.

We parametrize each boundary component $\Gamma_j$ by a $2\pi$-periodic complex function $\eta_j(t)$, $t\in J_j:=[0,2\pi]$, $j=1,2$. We assume that each of these functions $\eta_j(t)$ is twice continuously differentiable with $\eta'(t)\ne0$ (the presented method can be applied also if the curve $\Gamma_j$ has a finite number of corner points but no cusps~\cite{nmz}).
Then we define the total parameter domain $J$ as the disjoint union of the two intervals $J_1=[0,2\pi]$ and $J_2=[0,2\pi]$, i.e., $J =J_1\sqcup J_2=\cup_{j=1}^{2}\{(t,j)\;:\;t\in J_j\}$. The elements of the total parameter domain $J$ are ordered pairs $(t,j)$ where $t$ is a real number in $[0,2\pi]$ and the index $j$ is an integer indicating the interval $J_j$ containing $t$~\cite{Nas-ETNA}. Hence, the boundary $\Gamma$ can be parametrized by
\begin{equation}\label{e:eta-1}
\eta(t,j)=\eta_j(t), \quad t\in J_j,\quad j=1,2.
\end{equation}
For a given $t$, the index $k$ for which the interval $J_k$ contains $t$ will be always clear from the context, see e.g.,~\cite{Nas-Siam1,Nas-ETNA,NF13,NG18,Weg-Nas}. So the pair $(t,k)$ in the left-hand side of~(\ref{e:eta-1}) will be replaced by $t$ and a parametrization of the whole boundary $\Gamma$ can be defined on $J$ by
\begin{equation}\label{eq:eta}
\eta(t)= \left\{ \begin{array}{l@{\hspace{0.5cm}}l}
\eta_1(t),&t\in J_1,\\
\eta_2(t),&t\in J_2. \\
\end{array}
\right.
\end{equation}
We denote by $H$ the space of all functions of the form
\[
\rho(t)= \left\{ \begin{array}{l@{\hspace{0.5cm}}l}
\rho_1(t),&t\in J_1,\\
\rho_2(t),&t\in J_2, \\
\end{array}
\right.
\]
where $\rho_1(t)$ and $\rho_2(t)$ are $2\pi$-periodic H\"older continuous real functions on $J_1$ and $J_2$, respectively. 

Let $A$ be the complex function~\cite{Nas-ETNA}
\begin{equation}\label{eq:A}
A(t)= \left\{ \begin{array}{l@{\hspace{0.5cm}}l}
e^{\i(\frac{\pi}{2}-\theta(t))}(\eta(t)-\alpha),&{\rm if\;} G{\rm \; is\; bounded},\\
e^{\i(\frac{\pi}{2}-\theta(t))} ,&{\rm if\;} G {\rm\; is\; unbounded},\\
\end{array}
\right.
\end{equation}
where $\theta$ is a real function with constant value on each interval $J_j$, i.e., 
\[
\theta(t)=\theta_j \qquad {\rm for}\quad t\in J_j
\]
and $\theta_j$ is a real constant, $j=1,2$.
The generalized Neumann kernel $N(s,t)$ is defined for $(s,t)\in J\times J$ by
\begin{equation}\label{eq:N}
N(s,t) :=
\frac{1}{\pi}\Im\left(\frac{A(s)}{A(t)}\frac{\dot\eta(t)}{\eta(t)-\eta(s)}\right).
\end{equation}
The kernel $N(s,t)$ is continuous~\cite{Weg-Nas}. Hence, the integral operator $\bN$ defined by
\[
\bN\rho(s) := \int_J N(s,t) \rho(t) dt, \quad s\in J,
\]
is compact. The integral equation with the generalized Neumann kernel involves also the following kernel 
\begin{equation}\label{eq:M}
M(s,t) :=
\frac{1}{\pi}\Re\left(\frac{A(s)}{A(t)}\frac{\dot\eta(t)}{\eta(t)-\eta(s)}\right), \quad (s,t)\in J\times J,
\end{equation}
which has a singularity of cotangent type~\cite{Weg-Nas}. The integral operator $\bM$ defined on $H$ by
\[
\bM\rho(s) := \int_J  M(s,t) \rho(t) dt, \quad s\in J,
\]
is singular, but is bounded on $H$~\cite{Weg-Nas}. For more details, see~\cite{Weg-Nas}.

For the above function $A$ defined by~\eqref{eq:A}, the following integral equation
\begin{equation}\label{eq:ie}
(\bI-\bN)\rho=-\bM\gamma
\end{equation}
is uniquely solvable for any real function $\gamma$ in $H$~\cite{Nas-jmaa11}. Furthermore, if $\rho$ is the unique solution of the boundary integral equation~\eqref{eq:ie}, then the real function $h$ defined by
\begin{equation}\label{eq:h}
h=[\bM\rho-(\bI-\bN)\gamma]/2
\end{equation}
is a piecewise constant function on the boundary $\Gamma$, i.e.,
\[
h(t)=h_j \qquad {\rm for}\quad \eta(t)\in\Gamma_j
\]
where $h_j$ is a real constant, $j=1,2$~\cite{Nas-jmaa11}. Moreover, 
\begin{equation}\label{eq:f}
f(\eta(t))=\frac{\gamma(t)+h(t)+\i\rho(t)}{A(t)}, \quad \eta(t)\in\Gamma,
\end{equation}
are boundary values of an analytic function $f$ in the doubly connected domain $G$ such that $f(\infty)=0$ when $G$ is unbounded. For more details, see~\cite{Nas-jmaa11,Nas-ETNA} and the references cited therein.

\subsection{Numerical solution of the integral equation}

The MATLAB function \verb|fbie| in~\cite{Nas-ETNA} provides us with an efficient and fast method for solving the boundary integral equation~\eqref{eq:ie}. 
The function \verb|fbie| is based on discretizing the boundary integral equation~(\ref{eq:ie}) using the Nystr\"om method with the trapezoidal rule~\cite{Akt,Atk97,Tre-Trap}. This discretization leads to a non-symmetric linear system. 
Then, the MATLAB function $\mathtt{gmres}$ is used to solve the linear system. 
The matrix-vector multiplication in the GMRES method is computed using the MATLAB function $\mathtt{zfmm2dpart}$ in the toolbox $\mathtt{FMMLIB2D}$~\cite{Gre-Gim12}. The function \verb|fbie| provides us also with approximations to the piecewise constant function $h$ in~\eqref{eq:h}. 
The computational cost for the overall method is $O(n\log n)$ operations where $n$ (an even positive integer) is the number of nodes in each of the intervals $J_1$ and $J_2$.

For the accuracy of the obtained numerical results, it is known that the order of the convergence of the Nystr\"om method depends on the order of convergence of the used quadratic method~\cite{Atk97}. 
The quadratic method used in the function \verb|fbie| is the trapezoidal rule which gives surprisingly accurate numerical results for periodic functions~\cite{Atk97,Tre-Trap}. 
In view of of~\eqref{eq:N} and~\eqref{eq:M}, the smoothness of the two kernels $N(s,t)$ and $M(s,t)$ depends on the smoothness of the parametrization function $\eta(t)$. Similarly, in this paper, the function $\gamma$ on the right-hand side of the integral equation~\eqref{eq:ie} will be defined in terms of $\eta(t)$. 
Hence, the smoothness of the function $\gamma$ will depend also on the smoothness of the boundary $\Gamma$. 
Thus, the order of convergence of the trapezoidal rule depends on the smoothness of the boundary $\Gamma$ of the domain $G$. 
For domain with smooth boundaries, we use the trapezoidal rule with equidistant nodes.
The integrand in the integral equation~\eqref{eq:ie} will be $C^\infty$ smooth if the boundaries of the domains are $C^\infty$ smooth.
Hence the rate of convergence of the numerical method is $O(e^{-cn})$ with a positive constant $c$ (see~\cite[p.~223]{Kre14}). If the boundary is $C^{q+2}$ smooth $(q\ge0)$, then the rate of convergence of the numerical method is $O(1/n^q)$~\cite{Kre90}.
For domains with corners (excluding cusps), the derivatives of the solution $\rho(t)$ of the boundary integral equation~\eqref{eq:ie} have a discontinuity at the corner points. 
Thus, only poor convergence can be achieved if the trapezoidal rule with equidistant nodes is used. 
For such domains, accurate results can be obtained if we use the trapezoidal rule with a graded mesh~\cite{Kre90}. 
Such a graded mesh can be obtained by substituting a suitable new variable in the integral equation such that the discontinuity in the derivatives of $\rho(t)$ is removed~\cite{Kre90,LSN17}.

To use the MATLAB function \verb|fbie|, the vectors \verb|et|, \verb|etp|, \verb|A|, and \verb|gam| that contain the discretizations of the functions $\eta(t)$, $\eta'(t)$, $A(t)$, and $\gamma(t)$, respectively, will be stored in MATLAB. Then we call the function
\[
[\verb|rho|,\verb|h|] = \verb|fbie|(\verb|et|,\verb|etp|,\verb|A|,\verb|gam|,\verb|n|,\verb|iprec|,\verb|restart|,\verb|gmrestol|,\verb|maxit|)
\]
to compute the vectors \verb|rho| and \verb|h| which contain the discretizations of the solution of the integral equation $\rho(t)$ and the piecewise constant function $h(t)$, respectively. In the numerical experiments in this paper, we set the tolerances of the FMM and the GMRES method to be $0.5\times 10^{-14}$ and $10^{-14}$ by choosing $\mathtt{iprec}=5$ and $\mathtt{gmrestol}=10^{-14}$, respectively. 
We use the GMRES method without restart by choosing $\mathtt{restart}=[\,]$ and with the maximum number of iterations $\mathtt{maxit}=100$. 
The choice of the value of $n$ depends on the geometry of the domain $G$. 
If $G$ has a simple geometry and smooth boundary, we can obtain accurate numerical results by choosing moderate values of $n$. 
If $G$ has a complex geometry, for example if its boundary has corners or its boundary components are close to each other, it is required to choose a sufficiently large value of $n$ to obtain accurate results. 
For domains with corners, we choose $n$ as a multiple of the number of corners.
Once the discretizations of the two functions $\rho(t)$ and $h(t)$ are computed, we use \[
\verb|fet|=(\verb|gam|+\verb|h|+\verb|i|*\verb|rho|)./\verb|A|
\]
to find approximations to the boundary values of the function $f(z)$.
Then approximations to the values of the function $f(z)$ for any vector of points \verb|z| in $G$ can be obtained using the Cauchy integral formula.
Numerically we carry out this computation by applying the MATLAB function \verb|fcau|~\cite{Nas-ETNA} by calling
\[
\verb|fz| = \verb|fcau|(\verb|et|,\verb|etp|,\verb|fet|,\verb|z|)
\]
for bounded $G$ and by calling
\[
\verb|fz| = \verb|fcau|(\verb|et|,\verb|etp|,\verb|fet|,\verb|z|,\verb|n|,0)
\]
for unbounded $G$ (here $0=f(\infty)$).
For more details, we refer the reader to~\cite{Nas-ETNA}.

The computations presented in this paper were performed on ASUS Laptop with Intel(R) Core(TM) i7-8750H CPU @2.20GHz, 2208 Mhz, 6 Core(s), 12 Logical Processor(s), and 16GB RAM, using using MATLAB R2017a. The MATLAB tic toc commands were used to measure the computation times.

\subsection{Computing the conformal mapping for bounded domains}\label{sec:cm-b}

If the domain $G$ is bounded, then we can compute the conformal mapping $w=\Phi(z)$ from $G$ onto the annulus $\{w\in\C\;:\;q<|w|<1\}$ with the normalization
\[
\Phi(\alpha)>0
\]
as in the following theorem from~\cite{Nas-Siam1}. Here, $\alpha$ is a given auxiliary point in $G$.

\begin{thm}\label{thm:cm-an-b}
Let $\theta_1=\theta_2=\pi/2$, let the function $A$ be defined by~\eqref{eq:A}, and let the function $\gamma$ be defined by
\begin{equation}
\gamma(t)=-\log\left|\frac{\eta(t)-z_2}{\alpha-z_2}\right|, \quad t\in J.
\end{equation}
If $\rho$ is the unique solution of the boundary integral equation~\eqref{eq:ie} and the piecewise constant function $h$ is given by~\eqref{eq:h}, then the function $f$ with the boundary values~\eqref{eq:f} is analytic in the domain $G$, the conformal mapping $\Phi$ is
given by
\begin{equation}\label{eq:Phi-an-b}
\Phi(z)=e^{-h_1}\left(\frac{z-z_2}{\alpha-z_2}\right)e^{(z-\alpha)f(z)}, \quad z\in G\cup\Gamma,
\end{equation}
and the modulus $q$ is given by
\begin{equation}\label{eq:R-an-b}
q=e^{h_2-h_1}.
\end{equation}
\end{thm}

\subsection{Computing the conformal mapping for unbounded domains}\label{sec:cm-u}

For an unbounded domain $G$, the following theorem from~\cite{Nas-Siam1} provides us with a method to compute the conformal mapping $w=\Phi(z)$ from $G$ onto the annulus $\{w\in\C\;:\;q<|w|<1\}$ with the normalization
\[
\Phi(\infty)>0.
\]

\begin{thm}\label{thm:cm-an-u}
Let $\theta_1=\theta_2=\pi/2$, let the function $A$ be defined by~\eqref{eq:A}, and let the function $\gamma$ be defined by
\begin{equation}
\gamma(t)=-\log\left|\frac{\eta(t)-z_2}{\eta(t)-z_1}\right|, \quad t\in J.
\end{equation}
If $\rho$ is the unique solution of the boundary integral equation~\eqref{eq:ie} and the piecewise constant function $h$ is given by~\eqref{eq:h}, then the function $f$ with the boundary values~\eqref{eq:f} is analytic in the domain $G$ with $f(\infty)=0$, the conformal mapping $\Phi$ is
given by
\begin{equation}\label{eq:Phi-an-u}
\Phi(z)=e^{-h_1}\left(\frac{z-z_2}{z-z_1}\right)e^{f(z)}, \quad z\in G\cup\Gamma,
\end{equation}
and the modulus $q$ is given by
\begin{equation}\label{eq:R-an-u}
q=e^{h_2-h_1}.
\end{equation}
\end{thm}

\subsection{Computing the capacity of the doubly connected domain $G$}\label{sc:comp-cap}

Since the capacity is invariant under conformal mapping, we shall compute the capacity of the above doubly connected domain $G$ (for both cases, bounded and unbounded) by mapping $G$ onto the annulus $R=\{w\in\C\;:\;q<|w|<1\}$ using the method presented in the above two theorems. Then the capacity of $G$ is the same as the capacity of the annulus $R$ which is given by the formula
\begin{equation}\label{eq:cap-G-R}
  {\rm cap} (G)={\rm cap}(R)=\frac{2\pi}{{\rm log}(1/q)} \,.
\end{equation}

A MATLAB implementation of the above described method for computing the radius $q$ of the inner circle of the annulus $R=\{w\in\C\;:\;q<|w|<1\}$ and hence the capacity $\capa(G)$ for both bounded and unbounded doubly connected domains $G$ is given in the following function.
The actual values of the auxiliary points $\alpha$, $z_1$, and $z_2$ in~\eqref{eq:Phi-an-b} and~\eqref{eq:Phi-an-u} are not important provided that we choose these points to be sufficiently far away from the boundary of the domains $G$.

\begin{lstlisting}
function [q,cap] = annq (et,etp,n,zz,z2,type)
%
if type=='b' 
    alpha = zz; A = et-alpha; gam = -log(abs((et-z2)./(alpha-z2)));
elseif type=='u' 
    z1 = zz; A = ones(size(et));  gam = -log(abs((et-z2)./(et-z1)));
end
[~,h] =  fbie(et,etp,A,gam,n,5,[],1e-14,200);
q     =  exp(mean(h(n+1:2*n))-mean(h(1:n)));
cap   =  2*pi/log(1/q);
end
\end{lstlisting}

\section{Rings with  piecewise smooth boundaries}

The method described in the previous section will be used in this section to compute the capacity of several doubly connected domains $G$ with piecewise smooth boundaries. For the first two examples, the exact values of the capacity are known.

\subsection{Two confocal ellipses}

In this example, we consider the bounded doubly connected domain $G$ in the interior of the ellipse
\[
\eta_1(t)=\frac{1}{2}\left(r_1e^{\i t}+\frac{1}{r_1}e^{-\i t}\right), \quad 0\le t\le 2\pi,
\]
and in the exterior of the ellipse
\[
\eta_2(t)=\frac{1}{2}\left(r_2e^{-\i t}+\frac{1}{r_2}e^{\i t}\right), \quad 0\le t\le 2\pi,
\]
where $r_1>r_2>1$. The domain $G$ is the image of the ring $q=r_2/r_1<|\zeta|<1$ under the Joukowski map
\[
z=\Phi(\zeta)=\frac{1}{2}\left(r_1\zeta+\frac{1}{r_1\zeta}\right).
\]
Hence, the exact conformal capacity of $G$ is $\capa(G)=2\pi/\log(1/q)=2\pi/\log(r_1/r_2)$.

\begin{figure}[h] %
\centerline{\scalebox{0.4}{\includegraphics[trim=0 0 0 0,clip]{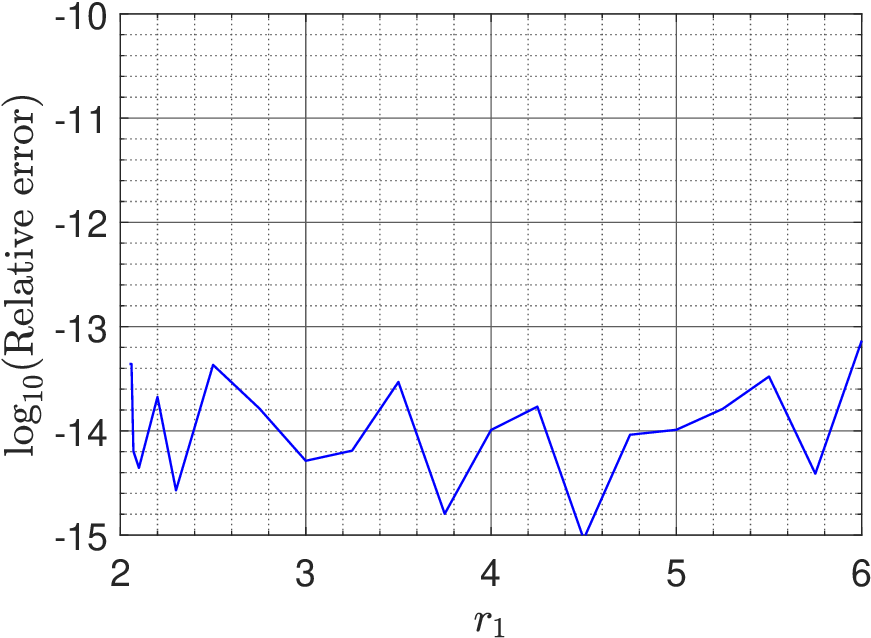}}
\hfill
\scalebox{0.4}{\includegraphics[trim=0 0 0 0,clip]{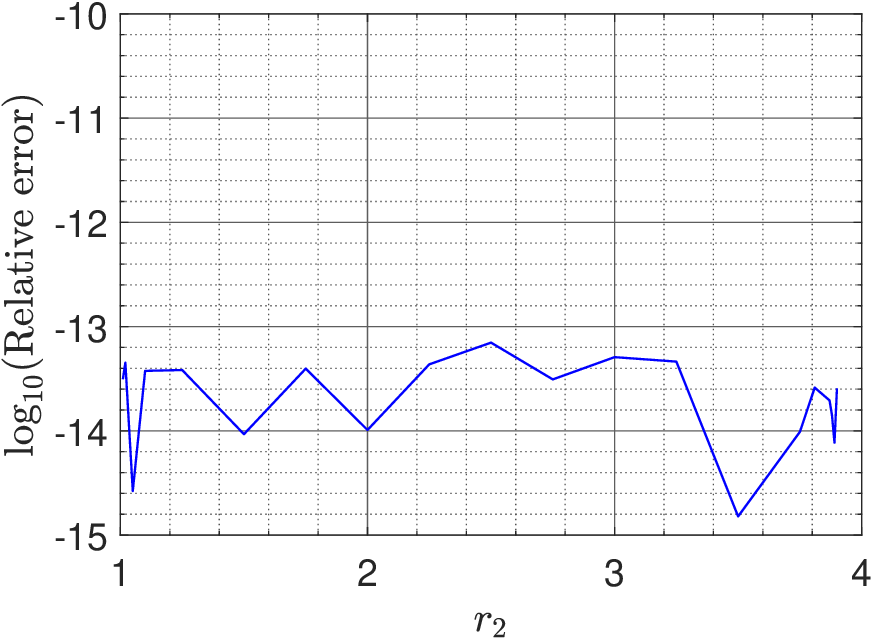}}}
\caption{Results for the two confocal ellipses domain: Relative errors in the computed conformal capacity for fixed $r_2=2$ (left) and fixed $r_1=4$ (right).}
\label{fig:two-ell}
\end{figure}

We use the MATLAB function \verb|annq| with $\alpha=((r_1+1/r_1)+(r_2+1/r_2))/4\in G$, $z_2=0$, and $n=2^{12}$ to calculate the capacity for several values of $r_1$ and $r_2$. First, we fixed $r_2=2$ and chose values of $r_1$ between $2.05$ and $6$. Then, we fixed $r_1=4$ and chose values of $r_2$ between $1.01$ and $3.9$. Figure~\ref{fig:two-ell} presents the relative errors in the computed values.

\subsection{Complete elliptic integrals}\label{mudef}

We recall the following facts about complete elliptic integrals and hypergeometric functions, needed for the sequel.
The {\it Gaussian hypergeometric function} is the analytic continuation of the series
\begin{equation} \label{eq:hypdef}
F(a,b;c;z) = {}_2 F_1(a,b;c;z) =
\sum_{n=0}^{\infty} \frac{(a,n)(b,n)}{(c,n)} \frac{z^n}{n!}\,,\:\: |z|<1 \,.
\end{equation}
 to the slit plane $\C \setminus [1,\infty)$ where $a,b,$ and $c$ are complex numbers with $c\neq 0,-1,-2, \ldots $.
Here $(a,n)$ is the {\it Appell symbol} or the {\it shifted factorial function}
$$
(a,n) = a(a+1)(a+2) \cdots (a+n-1)
$$
for $n \in \{1,2,3,\ldots\}$ and  $(a,0)=1$ for $a \neq 0$.
The
{\it complete elliptic integrals} of the first kind ${K}(r)$ and ${K}'(r)$ are defined by
\begin{equation} \label{defK}
{K}(r)=\frac{\pi}{2} F(1/2,1/2;1; r^2),
\qquad
{K}'(r)={K}(r'),\text{ and }r'= \sqrt{1-r^2}\,,
\end{equation}
and the {\it elliptic integrals} of the second kind ${E}(r)$ and ${E}'(r)$ are defined by 
\begin{equation} \label{defE}
{E}(r)=\frac{\pi}{2} F(1/2,-1/2;1; r^2),
\qquad
{E}'(r)={E}(r'),\text{ and }r'= \sqrt{1-r^2}\,.
\end{equation}
Then $K \colon (0,1)\to(\pi/2,\infty)$ is an increasing  homeomorphism and 
$E \colon (0,1)\to(1,\pi/2)$ is a decreasing  homeomorphism.
 The decreasing homeomorphism $\mu\colon (0,1)\to(0,\infty)$
is defined by
\begin{equation}
\label{mu}
\mu(r)= \frac{\pi}{2}\,\frac{{K}'(r)}{{K}(r)}\,.
\end{equation}

The basic properties of these functions can be found in \cite{hkv,avv,bf,olbc}. 
For example, it follows from~\cite[(5.2)]{avv} for $0<r<1$ that
\begin{equation}\label{eq:mu-r-double}
\mu(r)=2\mu\left(\frac{2\sqrt{r}}{1+r}\right)\,, \quad
\mu(r)=\frac{1}{2}\mu\left(\frac{1-r'}{1+r'}\right), \quad
\mu(r)\mu(r')=\frac{\pi^2}{4}.
\end{equation}
In the numerical calculations in this paper, we compute the values of $\mu(r)$ through~\eqref{mu} where the values of ${K}(r)$ and ${K}'(r)$ are computed by the MATLAB function \verb|ellipke|. Since $0<r<1$ and $r'=\sqrt{1-r^2}$, it readily follows that
\[
r<\frac{2\sqrt{r}}{1+r}<1\quad{\rm and}\quad 0<\frac{1-r'}{1+r'}<r.
\]
Thus, when $r$ is too close to $0$, we can use the first formula in~\eqref{eq:mu-r-double} to get accurate results with MATLAB function \verb|ellipke|. When $r$ is very close to $1$, we use the second formula in~\eqref{eq:mu-r-double}.

\subsection{Jacobi's inversion formula for $\mu$}\label{JacInv}

In his fundamental work on elliptic functions, C.G.J. Jacobi proved several dozens of formulas for these functions and related functions such as theta functions. Many of these formulas involved infinite products.
As pointed out in \cite[Thm 5.24(2)]{avv}, some of these formulas can
be rewritten so as to give formulas for $\mu^{-1}(y)\,.$ We give two
examples.
Jacobi's inversion formula for $\mu$  is  \cite[Thm 5.24(2)]{avv}
\[
\mu^{-1}(y)^2 =1- \prod\limits_{n=1}^\infty \left(
\frac{1-q^{2n-1}}{1+q^{2n-1}} \right)^8\,,\,  q= \exp(-2y)\, ,\,\, y > 0\,.
\]
Another example of Jacobi's work is the following formula for $\mu^{-1}(y)$ in terms of theta functions
\begin{equation}
 \label{invmutheta}
\mu^{-1}(y) = \left( \frac{\theta_2(0,q)}{\theta_3(0,q)}\right)^2\,,
\quad q= {\rm exp}(-2y) \,\,, y > 0\,,
\end{equation}
\begin{equation} \label{thetaDef}
\theta_2(0,q)=2 \sum_{n=0}^\infty q ^{(n+1/2)^2}\,, \quad
\theta_3(0,q)=1 + 2 \sum_{n=1}^\infty q ^{n^2}\,\,.
\end{equation}
Because these theta functions converge very fast in $[0, 0.95]\,,$ a few terms of series
expansion are enough to achieve numerical values correct up to 15 decimal places.
A Newton algorithm for computing $\mu^{-1}(y)$ was implemented
in \cite[pp. 92, 438]{avv}.

\subsection{Square in square}

In our second example, the domain $G$ is the difference of two concentric squares
 \[
 ((-2,2) \times (-2,2)) \setminus ((-2a,2a) \times (-2a,2a))
 \]
where $0<a<1$. 
The exact value of the capacity of this domain is \cite[pp. 103-104]{bow}
\begin{equation}\label{eq:cap-squ}
\capa(G)=\frac{4\pi}{\mu(r)}\,,
\end{equation}
where
\[
c=\frac{1-a}{1+a}, \quad u=\mu^{-1}\left(\frac{\pi c}{2}\right), \quad v=\mu^{-1}\left(\frac{\pi}{2c}\right), \quad
r=\left(\frac{u-v}{u+v}\right)^2\,.
\]
Then, by~\cite[Exercises 5.8(3)]{avv}, we have
\[
u^2+v^2=\left(\mu^{-1}\left(\frac{\pi c}{2}\right)\right)^2+\left(\mu^{-1}\left(\frac{\pi^2/4}{\pi c/2}\right)\right)^2=1
\]
and hence
\[
r=\frac{1-2uv}{1+2uv}\,.
\]
By~\cite[(5.2)]{avv}, we have 
\[
\mu(2uv)\mu(r)=\mu(2uv)\mu\left(\frac{1-2uv}{1+2uv}\right)=\frac{\pi^2}{2}\,.
\]
Thus, it follows from~\eqref{eq:cap-squ} that 
\begin{equation}\label{eq:cap-squ2}
\capa(G)=\frac{8}{\pi}\mu(2uv)\,,
\end{equation}

We use the MATLAB function \verb|annq| with $\alpha=1+a\in G$, $z_2=0$, and $n=2^{17}$ to calculate the capacity for several values of $a$ between $0.1$ and $0.9$. The obtained results are presented in Table~\ref{tab:HRV1}. Table~\ref{tab:HRV1} presents also the exact capacity and the numerical results computed in~\cite{hrv1} using an $hp$-FEM algorithm. We see from the results presented in the table that accurate results can be obtained using the presented method. The last column in Table~\ref{tab:HRV1} presents the CPU time (in seconds) for our method. The GMRES requires between $23$ to $25$ iterations only to converge. The obtained results using the presented method are not as accurate as the results obtained by the $hp$-FEM algorithm in~\cite{hrv1}. This is expected when we compare BIM and FEM for domains with corners.

\begin{table}[h]
\caption{The capacity values for the square in square domain.}
\label{tab:HRV1}%
\begin{tabular}{l|l|l|l|l}\hline
   $a$  &  Our Method     & \cite{hrv1}        & Exact value        &Time (sec)\\ \hline
0.1   & 2.83977741905231  & 2.83977741905223   & 2.83977741905224   & 6.6  \\
0.2   & 4.13448702423319  & 4.134487024234081  & 4.13448702423409   & 6.5  \\
0.3   & 5.63282800094106  & 5.632828000941654  & 5.63282800094165   & 6.5  \\
0.4   & 7.56153153980938  & 7.5615315398105745 & 7.56153153981058   & 7.1  \\
0.5   & 10.2340925693693  & 10.23409256936805  & 10.2340925693681   & 7.1  \\
0.6   & 14.2348796758222  & 14.234879675824363 & 14.2348796758244   & 6.6  \\
0.7   & 20.9015816764098  & 20.901581676413954 & 20.901581676414    & 6.4  \\
0.8   & 34.2349151987643  & 34.23491519877346  & 34.2349151987734   & 6.9  \\
0.9   & 74.2349151987441  & 74.23491519877882  & 74.2349151987788   & 6.9  \\
\hline %
\end{tabular}
\end{table}

\subsection{Polygon in polygon}

\begin{figure}[h] %
\centerline{\scalebox{0.4}{\includegraphics[trim=0 0 0 0cm,clip]{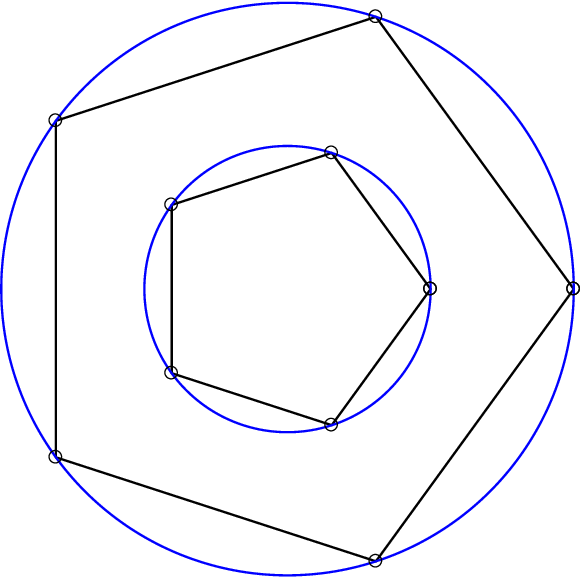}}}
\caption{The polygon in polygon domain $G$ for $m=5$ and $q=0.5$. The figure shows also the annulus $q<|z|<1$ where the vertices of the two polygons are on the circles $|z|=1$ and $|z|=q$.}
\label{fig:pg3}
\end{figure}

In the third example, we consider the doubly connected domain $G$ between two polygons. We assume that both polygons have $m$ vertices where $m\ge3$. We assume that the vertices of the external polygon are the roots of the unity and hence lie on the unit circle $|z|=1$. For the inner polygon, we assume that the vertices are the roots of the unity multiplied by $q=0.5$ and thus lie on the circle $|z|=q$ (see Figure~\ref{fig:pg3} for $m=5$).

The exact value of capacity of the domain is unknown (except for $m=4$ where the capacity can be computed as in the square in square example, which for $q=0.5$, is $10.2340925693681$).
We use the MATLAB function \verb|annq| with $\alpha=(1+q)/2\in G$, $z_2=0$, and $n=40320$ to calculate the capacity for several values of $m$. The computed capacity is presented in Table~\ref{tab:bsv1}. As we can see from the table, as $m$ increases, the capacity approaches the capacity of the annulus $q<|z|<1$ which is $2\pi/\log(1/q)$. For $q=0.5$, the capacity of the annulus is $9.064720283654388$. For some values of $m$, Table~\ref{tab:bsv1} presents also approximate values of the capacity from~\cite{bsv}. The last column in Table~\ref{tab:bsv1} presents the CPU time (in seconds) for our method.

\begin{table}[h]
\caption{The capacity values for the polygon in polygon domain.}
\label{tab:bsv1}%
\begin{tabular}{l|l|l|l}\hline
   $m$  &  Our Method         & \cite{bsv}  & Time (sec)\\ \hline
   3    &  12.4411574383      & 12.4412     & 4.0\\
   4    &  10.2340925693267   &             & 2.5\\
   5    &  9.62720096044514   & 9.6266      & 2.6\\
   7    &  9.25977557690559   & 9.2598      & 2.4\\
   9    &  9.15441235751744   & 9.1541      & 2.1\\
  15    &  9.08360686195382   &             & 1.8\\
	30    &  9.06705650051687   &             & 1.5 \\ \hline
\end{tabular}
\end{table}

\section{Complement of  two slits}
\label{sc:two-slits}

\begin{figure}[h] %
\centerline{
\hfill
\scalebox{0.35}{\includegraphics[trim=0 0 0 0cm,clip]{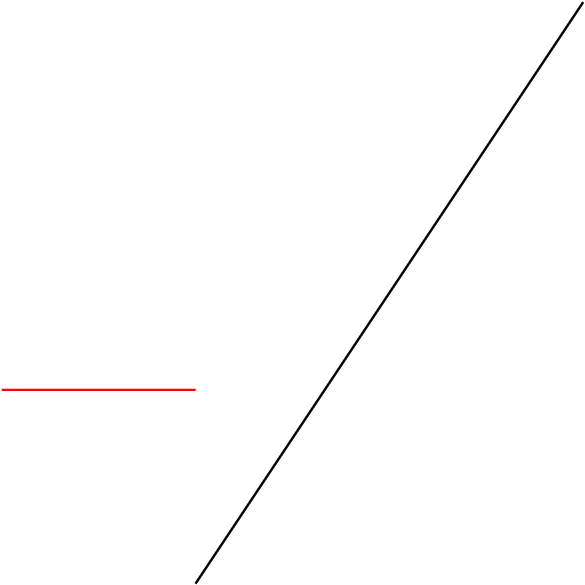}}
\hfill
\scalebox{0.35}{\includegraphics[trim=0 0 0 0cm,clip]{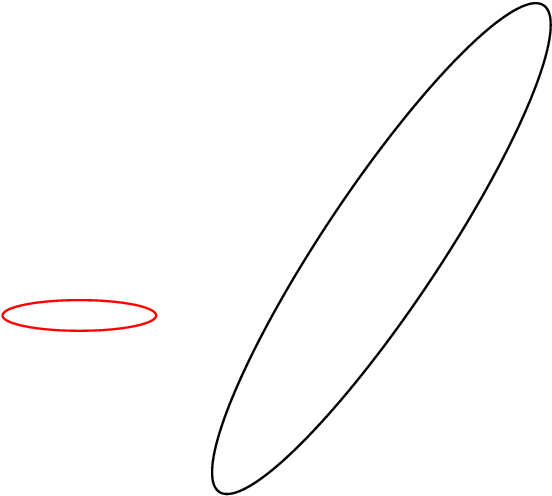}}
\hfill
}
\caption{The two segments domain $\Omega$ for $a=0$, $b=1$, $c=1-\i$ and $d=3+2\i$ (left) and the preimage domain $G$ bordered by ellipses (right).}
\label{fig:two-seg-fg}
\end{figure}

In this section, we consider a doubly connected domain $\Omega$ whose complementary components are the two non-intersecting segments $E=[a,b]$ and $F=[c,d]$ where $a$, $b$, $c$ and $d$ are complex numbers (see Figure~\ref{fig:two-seg-fg} (left) for $a=0$, $b=1$, $c=1-\i$ and $d=3+2\i$). 
Computing the capacity of such domain $\Omega$ has been considered recently in~\cite{dnv} using Weierstrass elliptic functions. Here, we shall compute the capacity of $\Omega$ using the method presented in Section~\ref{sec-cm}. 
However, a direct application of the method presented in Section~\ref{sec-cm} is not possible since the boundaries of $\Omega$ are not Jordan curves. 
So, we need to first map the given domain $\Omega$ onto a domain $G$ of the form considered in Section~\ref{sec-cm}. Up to the best of our knowledge, there is no analytic formula for a conformal mapping from the above doubly connected domain $\Omega$ onto a doubly connected domain $G$ bordered by smooth Jordan curves. So, we need to use numerical methods to find such an equivalent domain $G$. 
Such a conformally equivalent domain $G$ can be computed using the iterative method presented recently in~\cite{NG18}. 
The computed domain $G$ will be bordered by ellipses as in Figure~\ref{fig:two-seg-fg} (right). 
We refer the reader to~\cite{NG18} for details on this iterative numerical method. 
The MATLAB function \verb|annq| with $n=2^{11}$ is then used to compute the capacity of $G$, and hence the capacity of $\Omega$, for several values of $a$, $b$, $c$ and $d$, as in the following examples.

\subsection{Two segments on the real axis}\label{sec:2re-seg} 

\begin{figure}[h] %
\centerline{\scalebox{0.4}{\includegraphics[trim=0 0 0 0,clip]{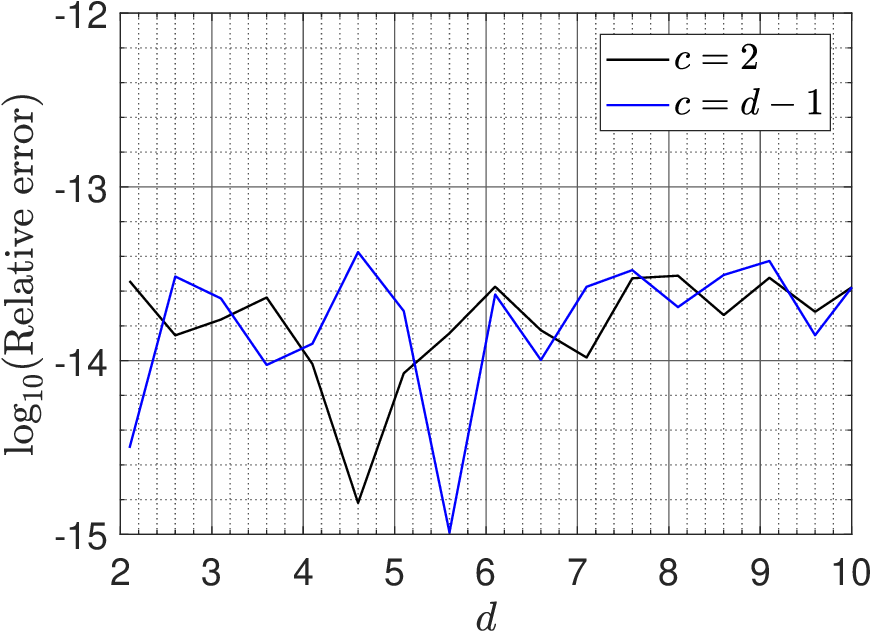}}}
\caption{Relative errors of the computed conformal capacity for ring domains with complementary components $[0,1] $ and $[c,d]\,,$ $2.1\le d\le10$.}
\label{fig:two-seg-err}
\end{figure}

When $E=[0,1]$ and $F=[c,d]$ with $d>c>1$, the exact capacity of $\Omega$ is known and is given by \cite[5.54 (1), 5.60(1)]{vu88}
\begin{equation}\label{eq:cap-two-seg}
\frac{\pi}{\mu\left(\sqrt{\frac{d-c}{c(d-1)}}\right)}\,.
\end{equation}
We tested our methods for several values of $c$ and $d$. First, we fixed $c=2$ and chose $d$ between $2.1$ and $10$. Then we fixed $c=d-1$ and chose $d$ between $2.1$ and $10$. For this case, the relative errors in the computed values are presented in Figure~\ref{fig:two-seg-err}.
As we can see from Figure~\ref{fig:two-seg-err}, the presented method gives accurate results with relative error around $10^{-14}$. Table~\ref{tab:two-seg-err} presents the approximate values of the capacity, the exact values of the capacity, and the total CPU time for several values of $c$ and $d$.

\begin{table}[h]
\caption{The values of $\capa([0,1],[c,d])$.}
\label{tab:two-seg-err}%
\begin{tabular}{l|l|l|l|l|l}\hline
  $c$    & $d$  & Computed value      & Exact value        & Relative Error      & Time (sec) \\ \hline %
$1.1$    & $2$  & $2.78768694945386$  & $2.7876869494539$  & $1.3\times10^{-14}$ & 3.8  \\
$1.1$    & $5$  & $3.11161184032646$  & $3.11161184032641$ & $1.7\times10^{-14}$ & 7.0  \\
$1.1$    & $10$ & $3.19100134481022$  & $3.19100134481039$ & $5.2\times10^{-14}$ & 10.0  \\
$2$      & $3$  & $1.56340192269607$  & $1.56340192269611$ & $2.7\times10^{-14}$ & 1.7  \\
$2$      & $5$  & $1.78056882835563$  & $1.78056882835559$ & $1.8\times10^{-14}$ & 2.5  \\
$2$      & $10$ & $1.9006702400055$   & $1.90067024000545$ & $2.5\times10^{-14}$ & 2.9  \\
\hline
\end{tabular}
\end{table}

\subsection{Two vertical segments}

The case $E=[a,\overline{a}]$ and $F=[c,\overline{c}]$, with $\Im a\ne0$ and $\Im c\ne0$, has been considered in~\cite[Figure E]{bbgghv}. We use our method to calculate the capacity for the same values of $a$ and $c$ that considered in~\cite[Table~3]{bbgghv}. A comparison of the results computed using our method vs the method presented in~\cite{bbgghv} is given in Table~\ref{tab:BBGGHV1} where the last column presents the CPU time for our method.

\begin{table}[h]
\caption{The values of $\capa([a,\overline{a}], [c,\overline{c}])$.}
\label{tab:BBGGHV1}%
\begin{tabular}{l|l|l|l|l}\hline
  $a$     & $c$     & Our Method          & \cite{bbgghv}      & Time (sec)   \\  \hline %
$0+\i$    & $5+2\i$ & 1.569943666568835   & 1.56994325474948999 & 3.2  \\
$0+2\i$   & $5+2\i$ & 1.873067768653831   & 1.87306699654806386 & 2.9  \\
$0+3\i$   & $5+2\i$ & 2.082038279851203   & 2.08203777712328096 & 3.8  \\
$0+4\i$   & $5+2\i$ & 2.232598863252026   & 2.23259828277206300 & 4.5  \\
$0+5\i$   & $5+2\i$ & 2.341589037102932   & 2.34158897620030515 & 5.0  \\
$0+3\i$   & $5+3\i$ & 2.352412309035929   & 2.35241226225174034 & 3.7  \\
\hline
\end{tabular}
\end{table}

\subsection{Two general segments}

\begin{figure}[h] %
\centerline{
\scalebox{0.45}{\includegraphics[trim=0 0 0 0,clip]{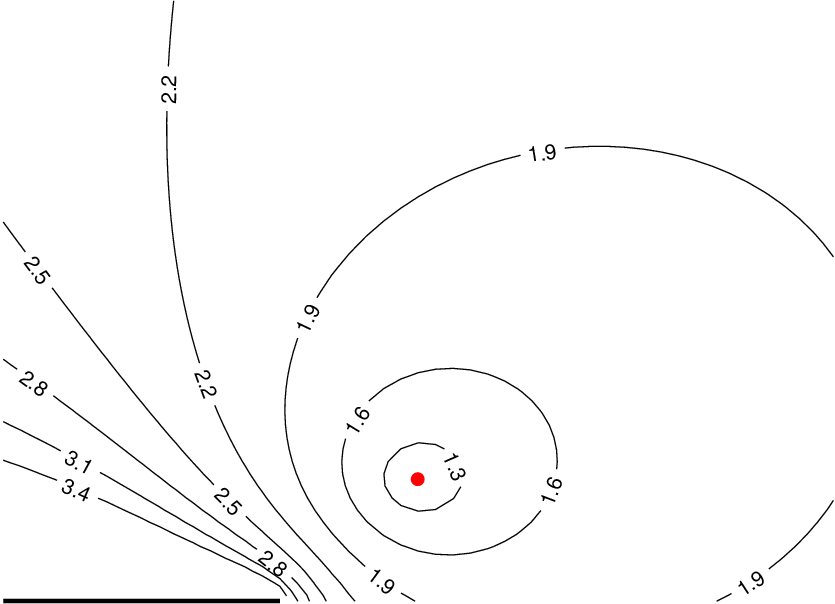}}
\hfill
\scalebox{0.45}{\includegraphics[trim=0 0 0 0,clip]{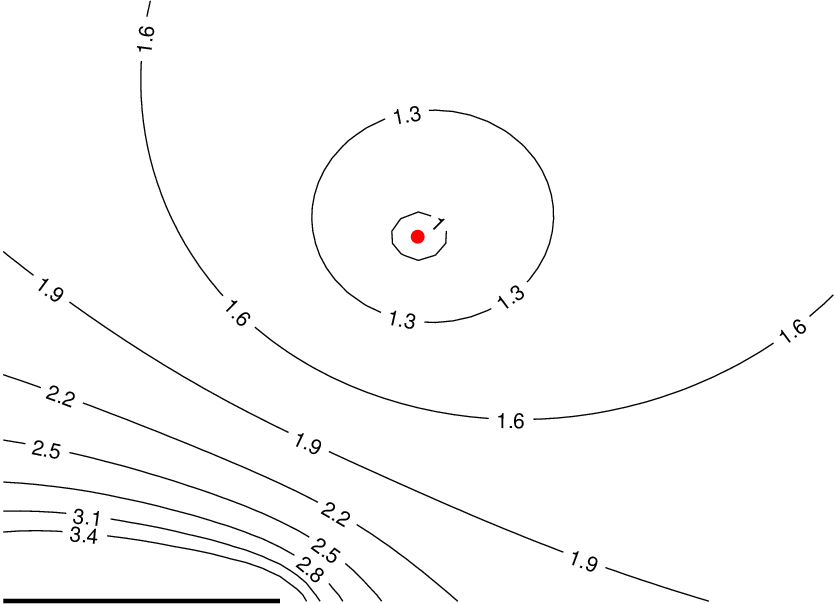}}
}
\caption{Results for the two segments domain: the contour lines of the function $u(x,y)=\capa([0,4],[6+2\i,x+\i y])$ (left) and $u(x,y)=\capa([0,4],[6+6\i,x+\i y])$ (right).}
\label{fig:two-seg-cant1}
\end{figure}

Finally, let
\[
f(a,b,c,d)=\capa([a,b],[c,d]),
\]
where $a$, $b$, $c$, and $d$ are complex numbers.
We fix $a=0$ and $b=4$.
Then, for a given point $z_1$ in the simply connected domain $\hat\Omega$ exterior to $[a,b]$, we define the function $u(x,y)$ by
\[
u(x,y)=f([0,4],[z_1,x+\i y]) =\capa([0,4],[z_1,x+\i y])\, .
\]
for $0<x<12$ and $0<y<10$ such that the segment $[z_1,x+\i y]$ is in $\hat\Omega$ with $x+\i y\ne z_1$. We plot the contour lines for the function $u(x,y)$ corresponding to several levels. The contour lines for $z_1=6+2\i$ and $z_1=6+6\i$ are shown in Figure~\ref{fig:two-seg-cant1}. Table~\ref{tab:cap-z1-z2-2s} presents the approximate values of $\capa([0,4],[z_1,z_2])$ for several values of $z_1$ and $z_2$.

\begin{figure}[h] %
\centerline{
\scalebox{0.45}{\includegraphics[trim=0 0 0 0,clip]{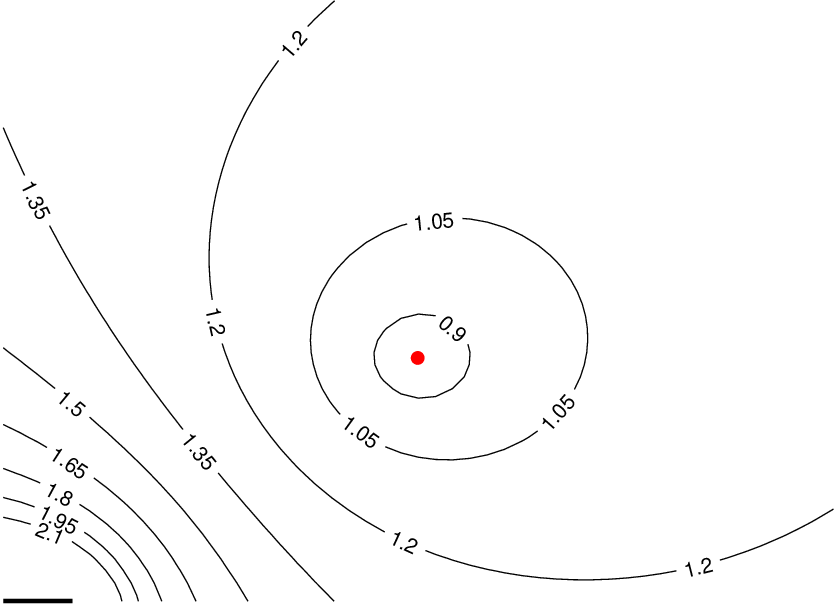}}
\hfill
\scalebox{0.45}{\includegraphics[trim=0 0 0 0,clip]{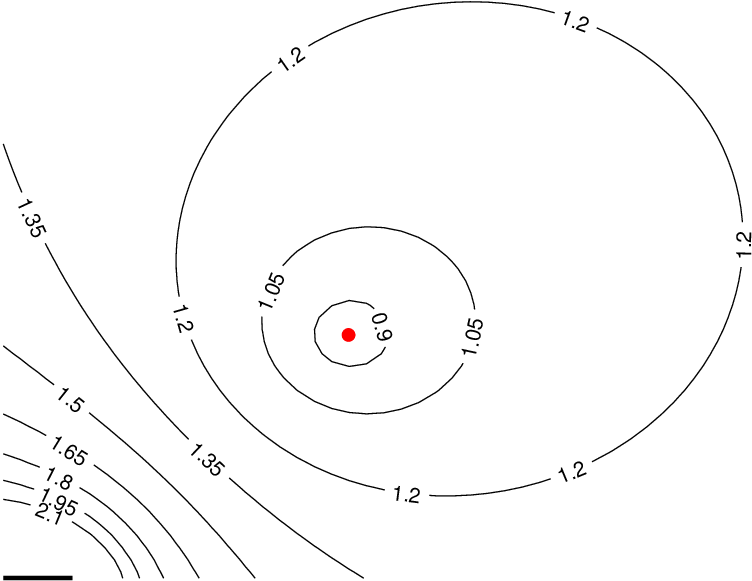}}
}
\caption{Results for the two segments domain: the contour lines of the function $u(x,y)=\capa([0,1],[6+4\i,x+\i y])$ (left) and $u(x,y)=\capa([0,1],[5+4\i,x+\i y])$ (right).}
\label{fig:two-seg-cant4}
\end{figure}

If the interval $[a,b]=[0,1]$ is considered instead of $[a,b]=[0,4]$, we obtain the results shown in Figure~\ref{fig:two-seg-cant4} for $z_1=6+4\i$ and $z_1=5+4\i$.

\begin{table}[h]
\caption{The values of $\capa([0,4],[z_1,z_2])$.}
\label{tab:cap-z1-z2-2s}%
\begin{tabular}{l|l|l|l}\hline
  $z_2\backslash z_1$   & $6+2\i$      & $6+4\i$      & $6+6\i$    \\  \hline %
 $1+\i$   & 4.437462457504561 & 3.780635179650131 & 3.564215562104226 \\
 $1+2\i$  & 3.317286587467568 & 2.860692915566007 & 2.711077789477010  \\
 $1+3\i$  & 2.846059598705353 & 2.436675855049381 & 2.295322432200487 \\
 $1+4\i$  & 2.604420470210280 & 2.202349785968325 & 2.046526840859631 \\
 $1+5\i$  & 2.470153941168786 & 2.066569200937597 & 1.886514461888595 \\
\hline
\end{tabular}
\end{table}

\section{Rings with a segment as a boundary component}

In this section, we compute the capacity of doubly connected domains $\Omega$ whose boundary components are  a slit and a piecewise smooth Jordan curve. Such domains cannot be mapped directly onto an annulus using the method presented in Section \ref{sec-cm}. To use the method presented in Section \ref{sec-cm}, we shall use first elementary mappings to map the domain $\Omega$ onto a domain $G$ of the types considered in Section \ref{sec-cm}. Then the domain $G$ is mapped onto an annulus $R=\{z\in\C\;:\;q<|z|<1\}$ and hence the capacity of $\Omega$ is $2\pi/\log(1/q)$. In this subsection we consider two examples where the exact value of the capacity for the first example is known.

\subsection{Segment and circle}\label{sec:cir}

\begin{figure}[ht] %
\centerline{
\hfill
\scalebox{0.35}{\includegraphics[trim=0 0 0 0cm,clip]{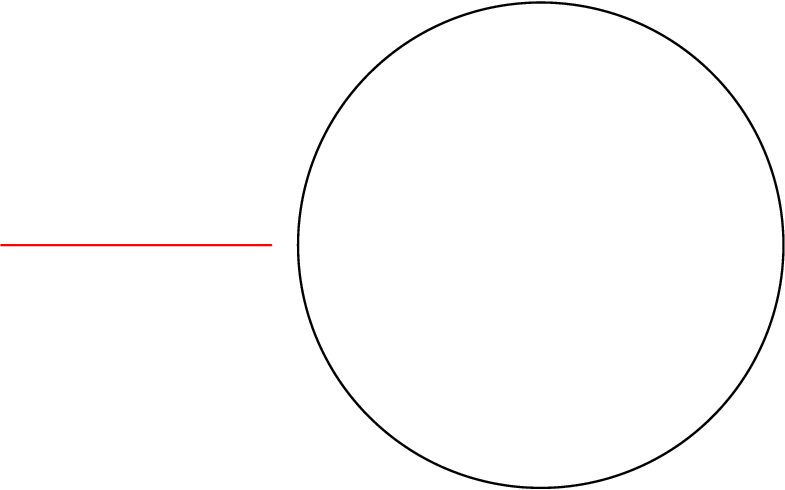}}
\hfill
\scalebox{0.35}{\includegraphics[trim=0 0 0 0cm,clip]{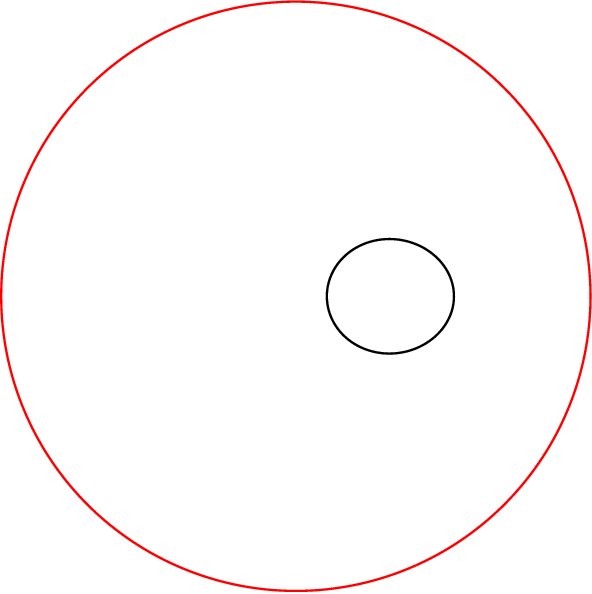}}
\hfill
}
\caption{The segment and circle domain $\Omega$ for $a=2$ and $r=0.9$ (left); and the image of this unbounded domain under the mapping $\zeta=\Psi^{-1}(z)$ (right).}
\label{fig:seg-cir-fg}
\end{figure}

First, we consider the doubly connected domain $\Omega$ in the exterior of the segment $[0,1]$ and the circle $|z-a|=r$ where $a$ is a real number with $a>1+r>1$ (see Figure~\ref{fig:seg-cir-fg} (left) for $r=0.9$ and $a=2$).
The exact conformal capacity for this domain $\Omega$ is known and given by~\cite[5.54(2)]{vu88}
\begin{equation}\label{eq:cap-disk-slit}
\capa(\Omega)=\frac{2\pi}{\mu(\tau)}, \quad \tau=\frac{r}{a^2-a-r^2}\,,
\end{equation}
where $\mu$ is given by~\eqref{mu}.

To apply our method presented in Section  \ref{sec-cm}, we shall use first elementary mappings to map the domain $\Omega$ onto a domain $G$ of the types considered in Section \ref{sec-cm}. It is known that the function
\[
z=\Psi(\zeta)=\frac{1}{4}\left(\zeta+\frac{1}{\zeta}\right)+\frac{1}{2}
\]
maps conformally the interior of the unit circle $|\zeta|=1$ onto the exterior of the segment $[0,1]$. Hence, its inverse function is given by
\begin{equation}\label{eq:Psi-inv}
\zeta=\Psi^{-1}(z)=\frac{1}{(2z-1)\left(1+\sqrt{1-\frac{1}{(2z-1)^2}}\right)},
\end{equation}
where we choose the branch for which $\sqrt{1}=1$. The function $\zeta=\Psi^{-1}(z)$ maps the segment $[0,1]$ onto the unit circle $|\zeta|=1$ and the exterior of the segment $[0,1]$ onto the interior of the unit circle $|\zeta|=1$. The function $\zeta=\Psi^{-1}(z)$ maps also the circle $|z-a|=r$ in the $z$-plane onto a smooth Jordan curve inside the unit circle $|\zeta|=1$. Consequently, the function $\zeta=\Psi^{-1}(z)$ maps the doubly connected domain $\Omega$ onto a bounded doubly connected domain $G$ of the form considered in Section \ref{sec-cm} (see Figure~\ref{fig:seg-cir-fg} (right)).

\begin{figure}[ht] %
\centerline{\scalebox{0.4}{\includegraphics[trim=0 0 0 0,clip]{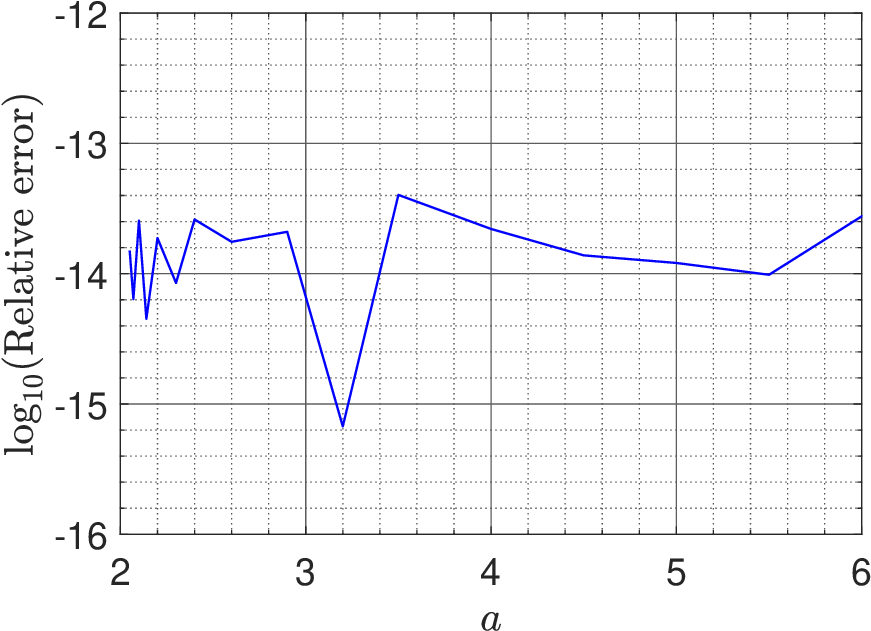}}
\hfill
\scalebox{0.4}{\includegraphics[trim=0 0 0 0,clip]{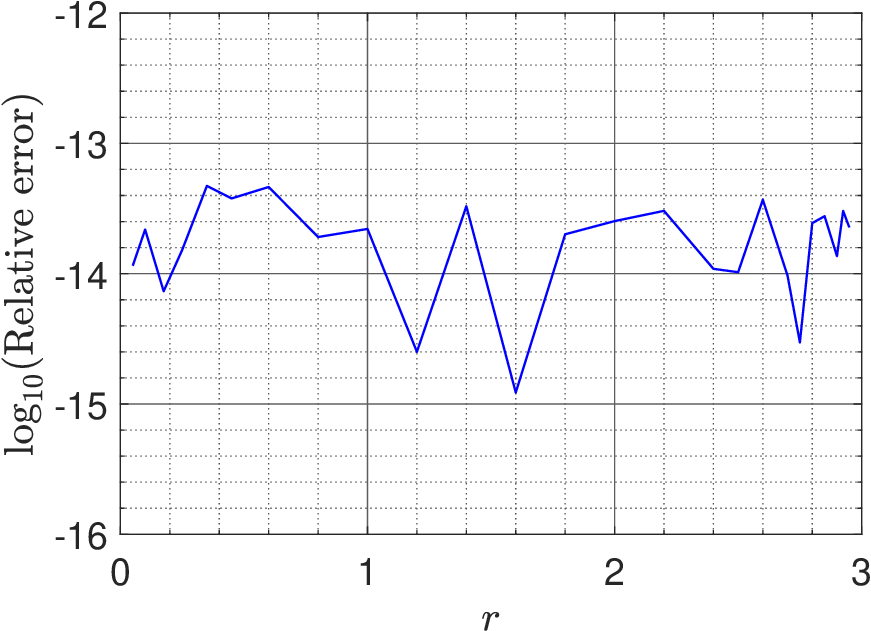}}}
\caption{Results for the segment and circle domain: Relative errors of the computed conformal capacity for fixed $r=1$ (left) and fixed $a=4$ (right).}
\label{fig:seg-cir}
\end{figure}

Then we use the MATLAB function \verb|annq| with $n=2^{11}$ to calculate approximate values for the capacity of $G$, and hence the capacity of $\Omega$, for several values of $a$ and $r$. First, we fixed $r=1$ and chose values of $a$ between $2.05$ and $6$. Then, we fixed $a=4$ and chose values of $r$ between $0.05$ and $2.95$.
Figure~\ref{fig:seg-cir} presents the relative errors in the calculated values. The exact values and the computed approximate values of the capacity are presented in Table~\ref{tab:seg-cir} for several values of $r$ and $a$.

\begin{table}[ht]
\caption{The capacity values for the segment and circle domain.}
\label{tab:seg-cir}%
\begin{tabular}{l|l|l|l|l|l}\hline
  $r$   & $a$   & Computed value  & Exact value & Relative Error & Time (sec)\\ \hline
  $0.1$ & $1.2$ & $2.89834979084902$  & $2.89834979084894$  & $2.7\times10^{-14}$ & 0.17\\
  $0.1$ & $2.2$ & $1.3496258349391$   & $1.34962583493908$  & $1.6\times10^{-14}$ & 0.15\\
  $0.1$ & $5.2$ & $0.927796431822476$ & $0.927796431822507$ & $3.3\times10^{-14}$ & 0.17\\
  $1.0$ & $2.1$ & $4.31652297947248$  & $4.31652297947259$  & $2.6\times10^{-14}$ & 0.19\\
  $3.0$ & $4.1$ & $4.6213142805315$   & $4.62131428053158$  & $1.8\times10^{-14}$ & 0.18\\
  $5.0$ & $6.1$ & $4.69478341049729$  & $4.69478341049717$  & $2.5\times10^{-14}$ & 0.19\\
\hline
\end{tabular}
\end{table}

\subsection{Segment and ellipse}\label{sec:ell}

\begin{figure}[ht] %
\centerline{
\scalebox{0.35}{\includegraphics[trim=0cm -3.4cm 0cm 0cm,clip]{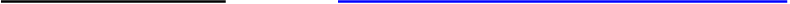}}
\hfill
\scalebox{0.35}{\includegraphics[trim=0cm -3.1cm 0cm 0cm,clip]{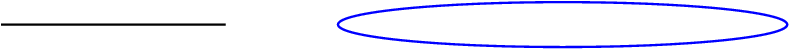}}
\hfill
\scalebox{0.35}{\includegraphics[trim=0cm 0cm 0cm 0cm,clip]{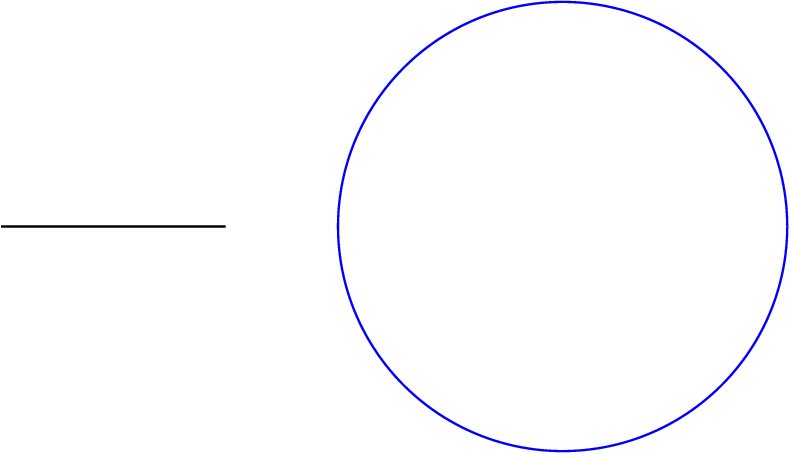}}
}
\caption{The domain $G_r$ with $c=1.5$, $d=3.5$  for $r=0$ (left), $r=0.1$ (center) and $r=b$ (right).}
\label{fig:slit-ellipse}
\end{figure}

In connection with the examples presented in Subsections~\ref{sec:2re-seg} and~\ref{sec:cir}, we consider the following example to show how the capacity of the domains changes when the geometry of the domains changes. 
Let $G_r$ be the doubly connected domain whose complementary components are the two non-intersecting closed sets $E=[0,1]$ and $F_r$ where $F_r$ is the closed set of points in the interior and on the boundary of the ellipse 
\[
\eta_r(t)=a+b\cos(t)-\i r\sin(t), \quad 0\le t\le2\pi, 
\]
where 
\[
a=\frac{1}{2}(d+c), \quad b=\frac{1}{2}(d-c), \quad 0<r\le b<a,
\]
and $1<c<d$ (see Figure~\ref{fig:slit-ellipse} (center)).

For $r=0$, $F_r$ reduced to the segment $F_0=[c,d]$ and hence $G_0$ is the doubly connected domain exterior to the two segments $E=[0,1]$ and $F_0=[c,d]$ as considered in Subsection~\ref{sec:2re-seg} (see Figure~\ref{fig:slit-ellipse} (left)). The exact value of $\capa(E,F_0)$ is given by~\eqref{eq:cap-two-seg}, i.e.,
\begin{equation}\label{eq:cap-s}
\capa(E,F_0)=\frac{\pi}{\mu\left(s\right)}\,, \qquad
s = \sqrt{\frac{d-c}{c(d-1)}},
\end{equation}
where $\mu$ is given by~\eqref{mu}. 
When $r=b$, $F_b$ is the closed disk $|z-a|\le b$ and the domain $G_b$ is then the doubly connected domain exterior to the segment $E=[0,1]$ and the closed disk $F_b=\{x\,|\,|z-a|\le b\}$ which was considered in Subsection~\ref{sec:cir} (see Figure~\ref{fig:slit-ellipse} (right)). The exact value of $\capa(E,F_b)$ is given by~\eqref{eq:cap-disk-slit}, i.e.,
\begin{equation}\label{eq:cap-disk-slit-c}
\capa(E,F_b)=\frac{2\pi}{\mu(s)}\,, \quad s=\frac{b}{a^2-a-b^2}\,.
\end{equation}

It is clear from the definition of the closed set $F_r$ that $F_0\subseteq F_r\subseteq F_b$ for $0\le r\le b$. As $r$ changes continuously from $0$ to $b$, the closed set $F_r$ changes continuously from the segment $F_0$ to the disk $F_b$. Here, we shall compute the exact value of the capacity $\capa(E,F_r)$ and show that $\capa(E,F_r)$ will change from $\capa(E,F_0)$ to $\capa(E,F_b)$ as $r$ changes from $0$ to $b$.

By the elementary mapping
\[
\zeta=\Psi_1(z)=\frac{z-a}{b},
\]
the unbounded domain $G$ is mapped conformally onto the unbounded domain $G_1$ exterior to the segment $[-a/b,-(a-1)/b]$ and the ellipse 
\[
\eta_1(t)=\cos(t)-\i (r/b)\sin(t), \quad 0\le t\le2\pi.
\]
We can easily show that the function
\[
\zeta=\Psi_2(\xi)=\xi+\frac{1-(r/b)^2}{4}\frac{1}{\xi}
\]
maps the domain exterior to the circle $|\xi|=(1+r/b)/2$ onto the domain exterior of the ellipse. 
Hence, the inverse mapping
\[
\xi=\Psi^{-1}_2(\zeta)=\zeta\left(\frac{1}{2}+\frac{1}{2}\sqrt{1-\frac{1-(r/b)^2}{\zeta^2}}\right),
\]
maps the domain $G_1$ onto the domain $G_2$ exterior to the circle $|\xi|=(1+r/b)/2$ and the segment $[c_1,d_1]$ where
\begin{equation}\label{eq:cap-tau-1}
c_1=-\frac{a+\sqrt{a^2-b^2+r^2}}{2b}, \quad
d_1=-\frac{a-1+\sqrt{(a-1)^2-b^2+r^2}}{2b},
\end{equation}
and the branch of the square root is chosen such that $\sqrt{1}=1$. Finally, the function 
\[
w=\Psi_3(\xi)=\frac{\xi-c_1}{d_1-c_1},
\]
maps the domain $G_2$ onto the domain $G_3$ exterior to the circle $|w-\hat a|=\hat r$ and the segment $[0,1]$ where
\begin{equation}\label{eq:cap-tau-2}
\hat a=-\frac{c_1}{d_1-c_1}, \quad
\hat r= \frac{b+r}{2b}\frac{1}{d_1-c_1}.
\end{equation}
Hence, the analytic value of $\capa(E,F_r)$ can be obtained since the exact value of conformal capacity of the domain $G_3$ is known~\cite[5.54(2)]{vu88},
\begin{equation}\label{eq:cap-tau2}
\capa(E,F_r)=\frac{2\pi}{\mu(\tau_r)}
\end{equation}
where
\begin{equation}\label{eq:cap-tau-3}
\tau_r=\frac{\hat r}{\hat a^2-\hat a-\hat r^2}\,.
\end{equation}
The value of $\tau_r$ can be obtained in terms of $c$, $d$ and $r$ as following
\[
\tau_r=\frac{2(d-c+2r)\left(1+\sqrt{dc+r^2}-\sqrt{dc-d-c+1+r^2}\right)}{\left(d+c+2\sqrt{dc+r^2}\right) \left(d+c-2+2\sqrt{dc-d-c+1+r^2}\right)-(d-c+2r)^2}.
\]

For $r=0$, the capacity given by~\eqref{eq:cap-tau2} becomes
\begin{equation}\label{eq:cap-tau0}
\frac{2\pi}{\mu(\tau_0)}
\end{equation}
where
\begin{eqnarray}
\nonumber
\tau_0&=&\frac{2(d-c)\left(1+\sqrt{dc}-\sqrt{dc-d-c+1}\right)}{\left(d+c+2\sqrt{dc}\right) \left(d+c-2+2\sqrt{dc-d-c+1}\right)-(d-c)^2}\\ \label{eq:tau-0}
      &=& \frac{\sqrt{d}-\sqrt{c}}{\sqrt{d}+\sqrt{c}} \,
			\frac{\sqrt{cd}-\sqrt{(c-1)(d-1)}+1}{\sqrt{cd}+\sqrt{(c-1)(d-1)}-1}
\end{eqnarray}
After tedious algebra, we find that $s$ in~\eqref{eq:cap-s} is related to $\tau_0$ in~\eqref{eq:tau-0} through
\begin{equation}\label{eq:s-tau0}
s=\frac{2\sqrt{\tau_0}}{1+\tau_0},
\end{equation}
which, in view of~\eqref{eq:mu-r-double}, implies that $\mu(\tau_0)=2\mu(s)$. Hence, 
\[
\capa(E,F_0)=\frac{\pi}{\mu(s)}=\frac{2\pi}{\mu(\tau_0)}
\]
and thus the capacity $\capa(E,F_r)$ given by~\eqref{eq:cap-tau2} reduced to the capacity $\capa(E,F_0)$ for $r=0$. Furthermore, when $r=b$, then it follows from~\eqref{eq:cap-tau-1} and~\eqref{eq:cap-tau-2} that $\hat a=a$ and $\hat r=r=b$. Hence, it follows from~\eqref{eq:cap-tau-3} that
\[
\tau_b=\frac{r}{a^2-a- r^2}\,,
\]
which implies that the capacity $\capa(E,F_r)$ given by~\eqref{eq:cap-tau2} reduced to the Formula~\eqref{eq:cap-disk-slit-c} for $r=b$. 

\begin{figure}[ht] %
\centerline{
\scalebox{0.4}{\includegraphics[trim=0cm 0cm 0cm 0cm,clip]{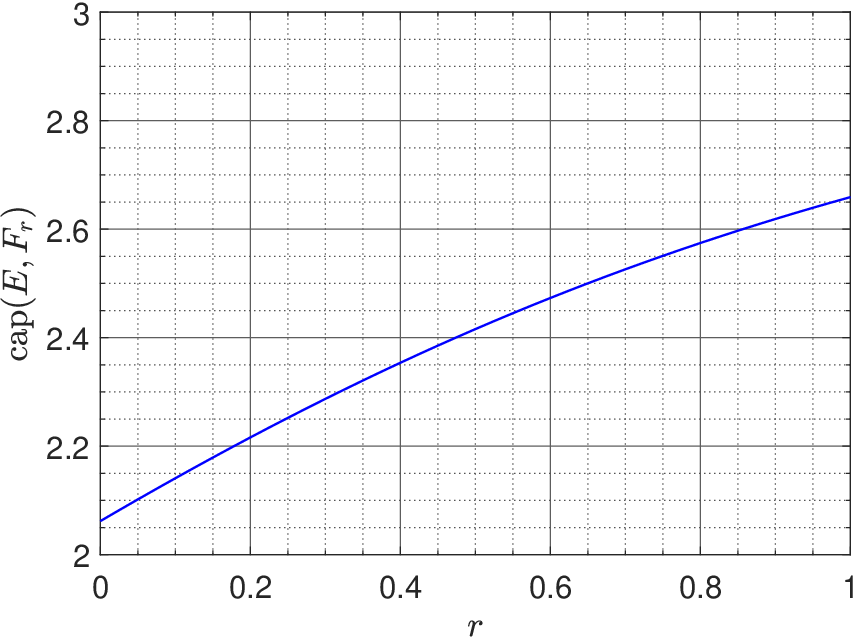}}
}
\caption{The values of $\capa(E,F_r)$ for $c=1.5$, $d=3.5$, $0\le r \le b$.}
\label{fig:slit-ellipse-r}
\end{figure}

The values of $\capa(E,F_r)$ for $c=1.5$, $d=3.5$, $0\le r \le b$ (where $b=(d-c)/2=1$) is given in Figure~\ref{fig:slit-ellipse-r}. As we can see from Figure~\ref{fig:slit-ellipse-r}, the capacity $\capa(E,F_r)$ changes continuously and rapidly increases from $\capa(E,F_0)$ to $\capa(E,F_b)$ as $r$ changes continuously from $0$ to $b$.

\subsection{Segment and polygon}

\begin{figure}[ht] %
\centerline{
\hfill
\scalebox{0.35}{\includegraphics[trim=+0cm -2cm 0cm 0cm,clip]{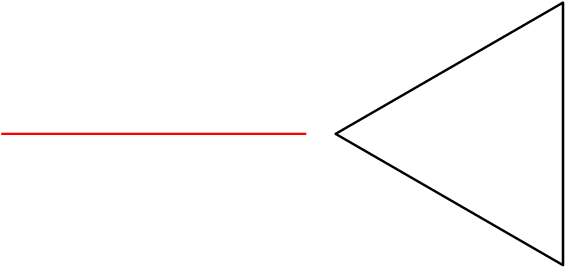}}
\hfill
\scalebox{0.35}{\includegraphics[trim=0 0 0 0cm,clip]{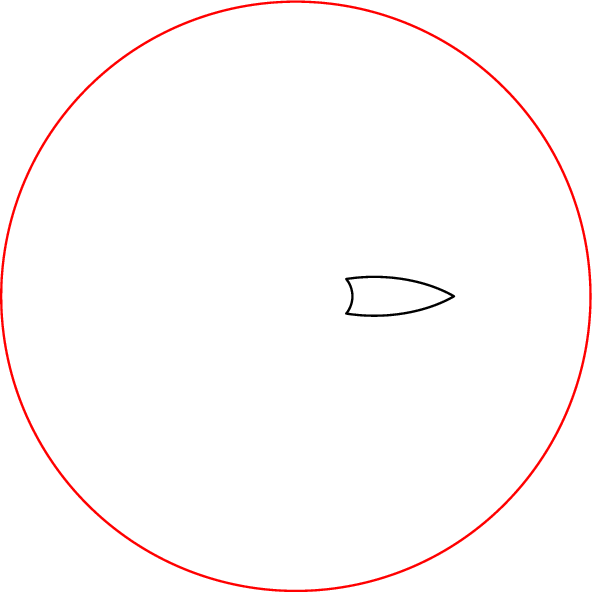}}
\hfill
}
\caption{The segment and polygon domain $\Omega$ for $a=1.6$ and $r=0.5$ (left); and the image of this unbounded domain under the mapping $\zeta=\Psi^{-1}(z)$ (right).}
\label{fig:seg-pol-fg}
\end{figure}

\begin{figure}[ht] %
\centerline{\scalebox{0.4}{\includegraphics[trim=0 0 0 0,clip]{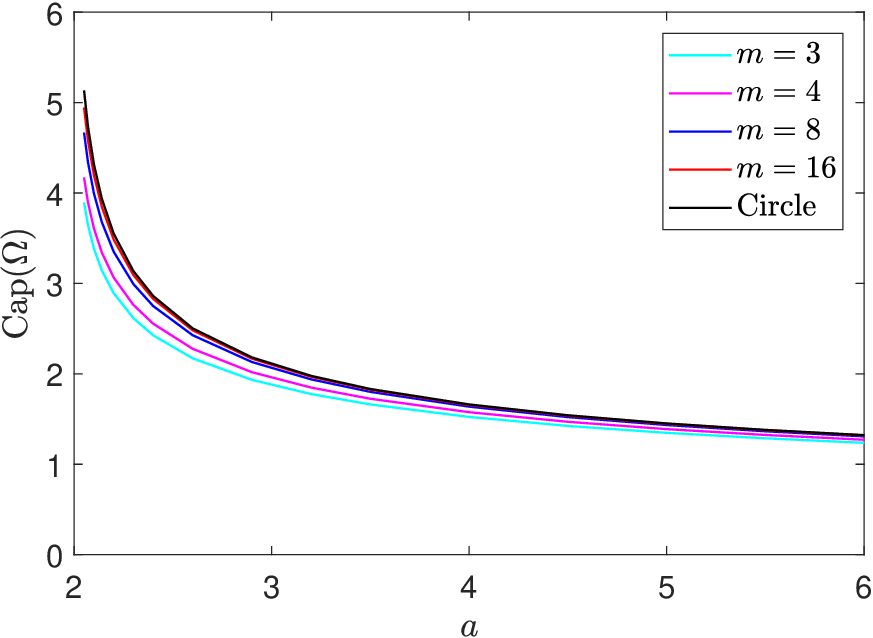}}
\hfill
\scalebox{0.4}{\includegraphics[trim=0 0 0 0,clip]{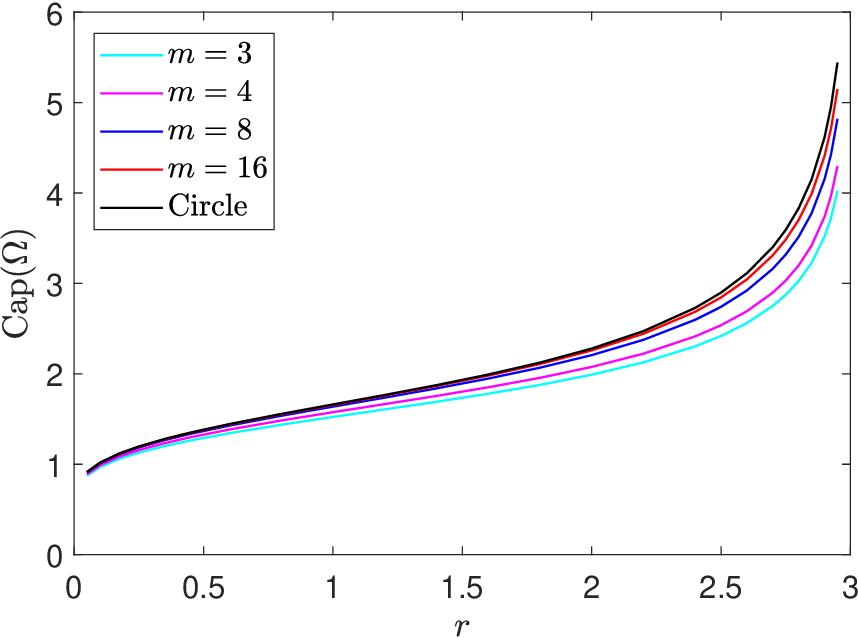}}}
\caption{Results for the segment and polygon domain: The computed conformal capacity for fixed $r=1$ (left) and fixed $a=4$ (right).}
\label{fig:seg-pol}
\end{figure}

In this example, we consider the doubly connected domain $\Omega$ in the exterior of the segment $[0,1]$ and a polygon with $m$ vertices where $m\ge3$. We assume that the vertices of the polygon are given by
\[
v_k=a-re^{\frac{-2k\pi\i}{m}}, \quad k=0,1,2,\ldots,m-1
\]
(see Figure~\ref{fig:seg-pol-fg} (left) for $a=1.6$, $r=0.5$, and $m=3$).
For this example, the exact value of the conformal capacity is unknown. To use the method described in Section~\ref{sec-cm}, we first use the mapping function $\zeta=\Psi^{-1}(z)$ in Subsection~\ref{sec:cir} to map the doubly connected domain $\Omega$ onto a bounded doubly connected domain $G$ of the form we considered in Section \ref{sec-cm} (see Figure~\ref{fig:seg-pol-fg} (right)). 
Then, for the new domain $G$, the MATLAB function \verb|annq| is used with $n=15\times2^9$ to calculate approximate values for the capacity of $\Omega$ for several values of $m$, $a$ and $r$. First, we fixed $r=1$ and chose values of $a$ between $2.05$ and $6$. The computed capacity for $m=3,4,8,16$ are presented in Figure~\ref{fig:seg-pol} (left). Then, we fixed $a=4$ and chose values of $r$ between $0.05$ and $2.95$. The computed capacity for $m=3,4,8,16$ are presented in Figure~\ref{fig:seg-pol} (right). Figure~\ref{fig:seg-pol} presents also the capacity for the segment with circle domain in the previous examples for the same values of $a$ and $r$. Table~\ref{tab:seg-cir} presents the calculated values of the capacity for the segment with circle domain and for the segment with polygon domain for several values of $a$, $r$, and $m$. As we can see from the results presented in the table, the capacity of the segment and polygon domain approaches the the capacity of the segment and circle domain as the number of vertices $m$ increases.

\begin{table}[ht]
\caption{The capacity values for the segment with polygon domain.}
\label{tab:seg-cir}
\begin{tabular}{l|l|l|l|l}\hline
 $r$   &$a$    &$m$   &Capacity (segment and polygon) &Capacity (segment and circle)   \\ \hline
 $1$   &$2.1$  &$3$   &$3.385465691885468$            &$4.31652297947259$  \\
       &       &$8$   &$3.996010644504850$            &   \\
       &       &$16$  &$4.198837938505387$            &   \\
       &       &$128$ &$4.314154067689326$            &   \\ \hline
 $0.5$ &$4$    &$3$   &$1.291427925789600$            &$1.38309579015095$  \\
       &       &$8$   &$1.368162812590014$            &   \\
       &       &$16$  &$1.379193284259540$            &   \\
       &       &$128$ &$1.383032359435526$            &   \\ \hline
 $2.5$ &$10$   &$3$   &$1.199970598794575$            &$1.28290663972126$  \\
       &       &$8$   &$1.268817744415183$            &   \\
       &       &$16$  &$1.279211726247828$            &   \\
       &       &$128$ &$1.282846509334037$            &   \\
	\hline
\end{tabular}
\end{table}

\section{The upper half-plane with a slit}
\label{sec:half-plane}

\begin{figure}[ht] %
\centerline{
\hfill
\scalebox{0.35}{\includegraphics[trim=0 -2.0cm 0 0cm,clip]{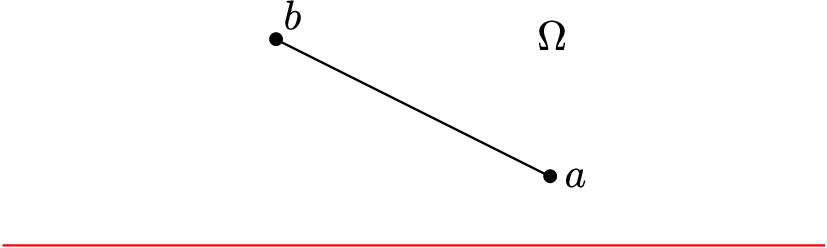}}
\hfill
\scalebox{0.35}{\includegraphics[trim=0 0 0 0cm,clip]{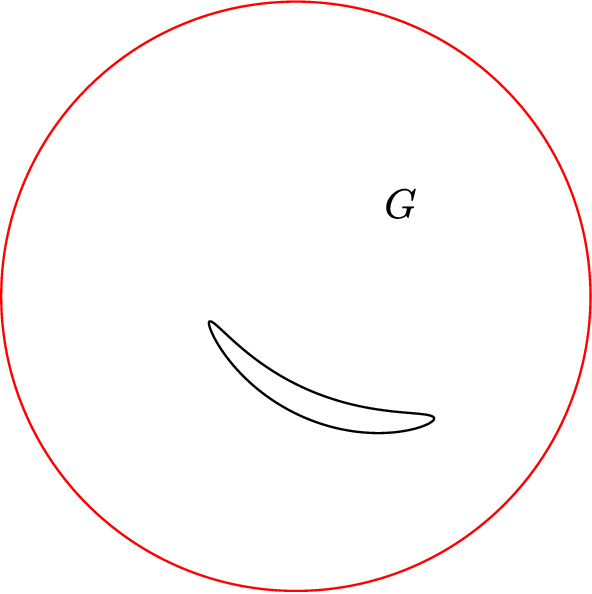}}
\hfill
}
\caption{The half-plane with a segment domain $\Omega$ for $a=1+0.5\i$ and $b=-1+1.5\i$ (left) and the preimage domain $G$ bordered by smooth Jordan curves (right).}
\label{fig:half-seg-fg}
\end{figure}

In this section, we consider the doubly connected domain $\Omega = \mathbb{H}^2 \setminus [a,b]$
where $\mathbb{H}^2$ is the upper half-plane $\{ z \in \mathbb{C}: \Im(z)>0\}$ and, $a$ and $b$ are two complex numbers in $\mathbb{H}^2$  (see Figure~\ref{fig:half-seg-fg} (left)).
For such domains $\Omega$, we cannot directly apply the method described in Section~\ref{sec-cm}. So, we first map the domain $\Omega$ onto a domain $G$ of the forms considered in Section \ref{sec-cm}. 
Since there is no exact conformal mapping from a domain such as $\Omega$ onto a doubly connected domain $G$ bordered by smooth Jordan curves, we find such an equivalent domain $G$ using numerical methods.
In this paper, we compute such a domain $G$ using the iterative numerical method presented in~\cite{NG18} (see Figure~\ref{fig:half-seg-fg}). We will omit the details here about the iterative method and refer the reader to~\cite{NG18}.
Then, we compute the capacity of the the given domain $\Omega$ by applying the MATLAB function \verb|annq| with $n=2^{11}$ to the new domain $G$.

\begin{figure}[ht] %
\centerline{
\scalebox{0.4}{\includegraphics[trim=0 0 0 0,clip]{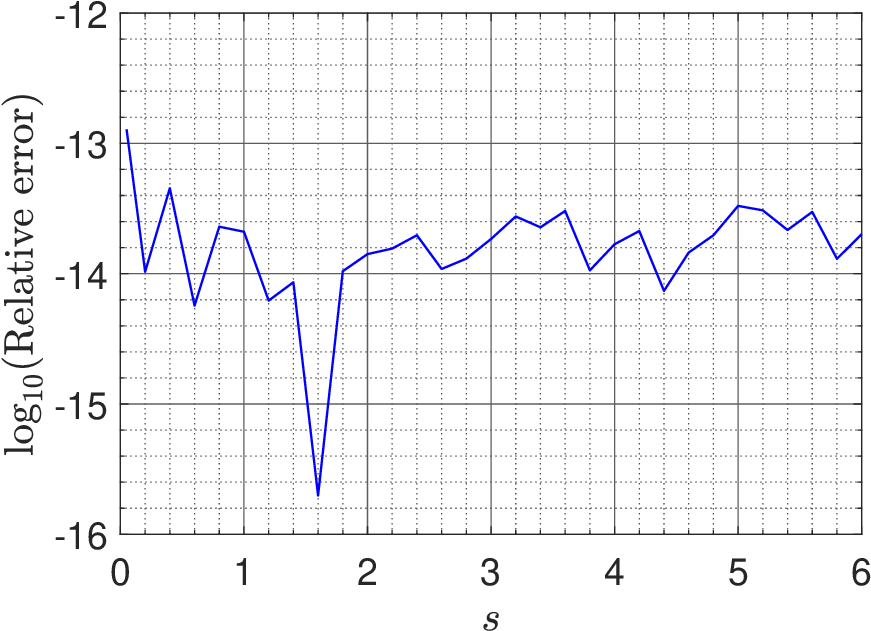}}
\hfill
\scalebox{0.4}{\includegraphics[trim=0 0 0 0,clip]{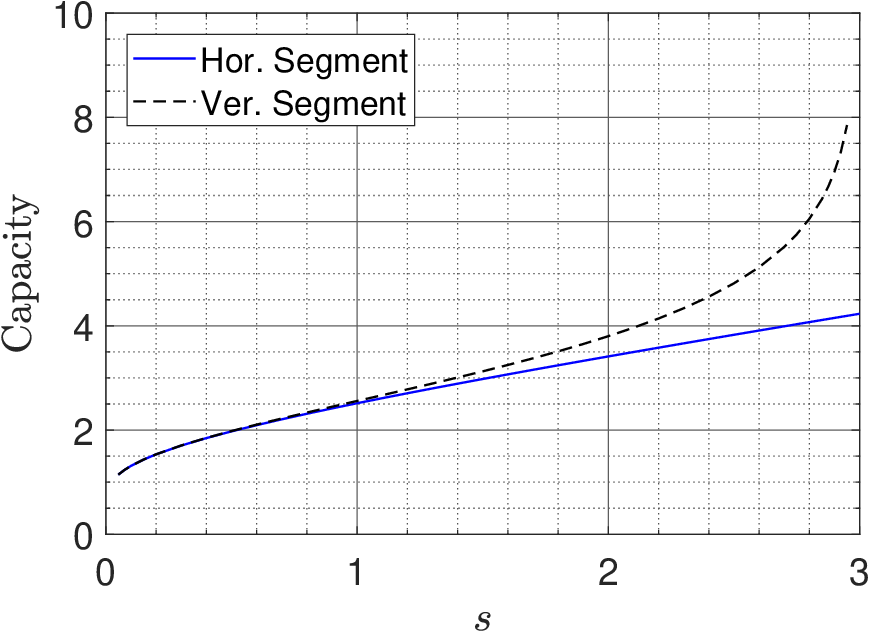}}
}
\caption{Results for the half-plane with a segment domain: Relative errors of the computed conformal capacity for the segment $F=[s\i,(1+s)\i]$ for $0.05\le s\le 6$ (left) and the computed capacities (right).}
\label{fig:half-seg-err}
\end{figure}

For the segment $F=[s\i,r\i]$ where $r>s>0$ are real numbers, the exact capacity of $\Omega$ is known and is given by \cite[(5.56), Theorem 8.6 (1)]{vu88}
\begin{equation}\label{eq:cap-half-seg}
\frac{2\pi}{\mu\left(\tanh\left(\frac{1}{2}\log\frac{r}{s}\right)\right)}\,.
\end{equation}
We tested our methods for several values of $s$ and $r$. First, we chose the vertical segment $F=[s\i,(1+s)\i]$, i.e., $a=s\i$ and $b=(1+s)\i$, for $0.05\le s\le 6$. For this case, the relative errors in the calculated values of the capacity are presented in Figure~\ref{fig:half-seg-err} (left). We see from Figure~\ref{fig:half-seg-err} (left) that the proposed method gives accurate results with relative error around $10^{-14}$.
The calculated and the exact values of the capacity as well as the total CPU time for several values of $s$ and $r$ are presented in Table~\ref{tab:half-seg-err}.

\begin{table}[ht]
\caption{The values of $\capa(\Omega)$ for $\Omega = \mathbb{H}^2 \setminus [s\i,r\i]$.}
\label{tab:half-seg-err}%
\begin{tabular}{l|l|l|l|l|l}\hline
  $s$    & $r$  & Computed value     & Exact value    & Relative Error      & Time (sec) \\ \hline %
$0.1$    & $1$  & $4.69363108974789$  & $4.6936310897475$  & $8.2\times10^{-14}$ & 2.9  \\
$0.1$    & $5$  & $6.74589984699685$  & $6.74589984699653$ & $4.8\times10^{-14}$ & 5.3  \\
$0.1$    & $10$ & $7.62853775997519$  & $7.62853775997481$ & $5.0\times10^{-14}$ & 6.6  \\
$1$      & $2$  & $2.55852314234207$  & $2.55852314234201$ & $2.1\times10^{-14}$ & 2.2  \\
$1$      & $5$  & $3.80134048001095$  & $3.80134048001091$ & $1.7\times10^{-14}$ & 2.6  \\
$1$      & $10$ & $4.6936310897476$   & $4.6936310897475$  & $2.1\times10^{-14}$ & 3.1  \\
\hline
\end{tabular}
\end{table}

We also compute the values of the capacity for the vertical segment $F=[(3-s)\i,(3+s)\i]$ for $0.05\le s\le2.95$ and for the horizontal segment $F=[-s+3\i,s+3\i]$ for $0.05\le s\le 3$. Both segments pass through the point $3\i$ and have the length $2s$. The results are presented in Figure~\ref{fig:half-seg-err} (right). Figure~\ref{fig:half-seg-err} (right) shows that the capacity increases as the length of the segment increases. For vertical segment, the capacity increases more rapidly when the segment becomes close to the real line.

\begin{figure}[ht] %
\centerline{
\scalebox{0.45}{\includegraphics[trim=0 0 0 0,clip]{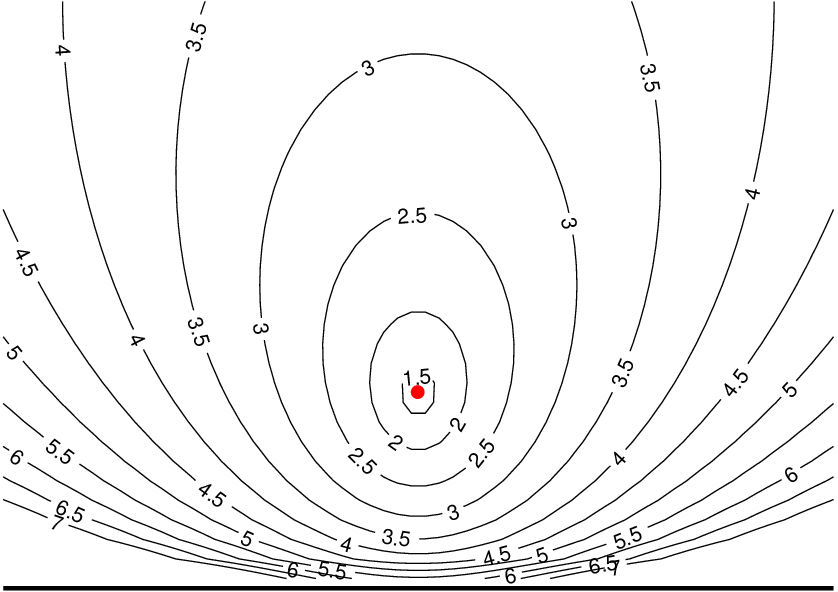}}
\hfill
\scalebox{0.45}{\includegraphics[trim=0 0 0 0,clip]{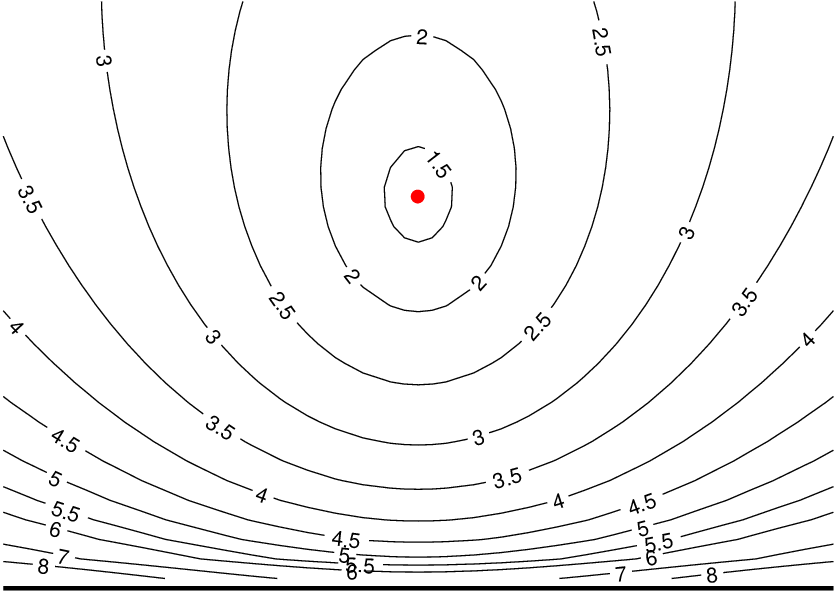}}
}
\caption{Results for the half-plane with a segment domain: the contour lines of the function $u(x,y)=\capa(\mathbb{H}^2\backslash[\i,x+\i y])$ (left) and $u(x,y)=\capa(\mathbb{H}^2\backslash[2\i,x+\i y])$ (right).}
\label{fig:half-seg-cont}
\end{figure}

Finally, for a given point $z_1$ in $\mathbb{H}^2$, we define the function $u(x,y)$ by
\[
u(x,y)=\capa(\mathbb{H}^2\backslash[z_1,x+\i y])
\]
for $-3<x<3$ and $0<y<3$ such that $x+\i y\ne z_1$. We plot the contour lines for the function $u(x,y)$ corresponding to several levels. The contour lines are shown in Figure~\ref{fig:half-seg-cont} for $z_1=\i$ and $z_1=2\i$.

\section{Domains exterior to thin rectangles}

\subsection{Two rectangles}\label{sec:two-rec}

\begin{figure}[ht] %
\centerline{
\scalebox{0.4}{\includegraphics[trim=0 0 0 0cm,clip]{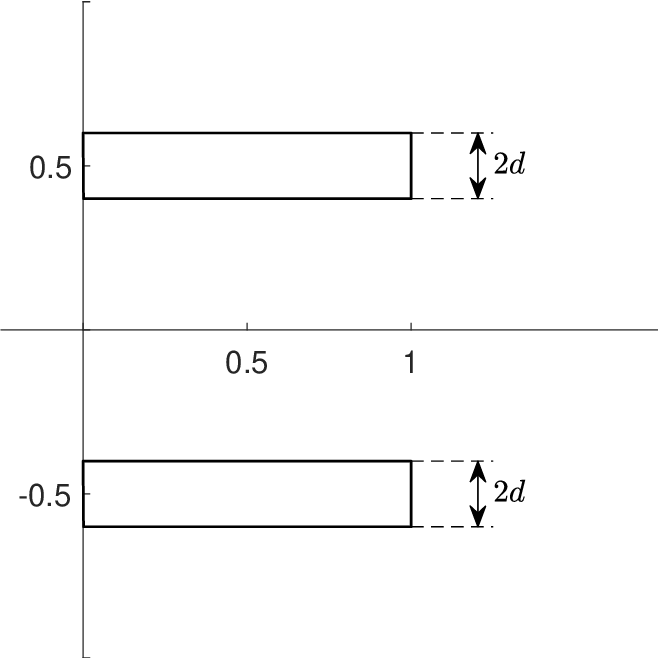}}
}
\caption{The domain $G$ in the exterior of the rectangular closed sets for $d=0.1$.}
\label{fig:two-rec}
\end{figure}

We consider in this section the doubly connected domain $G$ exterior to the rectangular closed sets
 \[
 [0,1] \times [0.5-d,0.5+d] \quad{\rm and}\quad [0,1] \times [-0.5-d,-0.5+d]
 \]
where $0<d<0.5$ (see the Figure~\ref{fig:two-rec}). 
We use the MATLAB function \verb|annq| presented in Subsection~\ref{sc:comp-cap} with $n=2^{15}$ to compute the capacity of $G$ for several values of $d$. 
When $d=0$, the two rectangles reduced to the two slits $[\i/2,1+\i/2]$ and $[-\i/2,1-\i/2]$. For these two slits, we can use the numerical method presented in Section~\ref{sc:two-slits} to compute the capacity of the domain in the exterior to these two slits. The obtained results are presented at the bottom of Table~\ref{tab:two-rec}.

\begin{table}[ht]
\caption{The capacity values for the domain exterior to the two rectangles.}
\label{tab:two-rec}%
\begin{tabular}{l|l|l}\hline
 $d$  & Capacity          & Time (sec)\\ \hline
0.4   & 7.55672805385065  & 2.1  \\
0.3		& 4.55284511607753	& 2.1  \\
0.2   & 3.3856923786737   & 2.0  \\
0.1   & 2.68688786213937  & 2.0  \\
0.05  & 2.40554719800866  & 2.1  \\
0.02  & 2.24063059387802  & 2.4  \\
0.01  & 2.18262548680027  & 2.8  \\
0.005 & 2.15161636330889  & 4.9  \\
0     & 2.11577897412447  & 25.0  \\
\hline
\end{tabular}
\end{table}

Let $R_1$ be the unbounded doubly connected domain exterior to the two slits $[\i/2,1+\i/2]$ and $[-\i/2,1-\i/2]$ (corresponding to $d=0$). The exact value of the capacity of $R_1$ can be computed. For $0<k<1$, consider the unbounded doubly connected domain $R_2$ exterior to the two slits $[-1/k,-1]$ and $[1,1/k]$. Then the M\"obius transform
\[
\Psi(z)=\frac{2k}{k-1}\frac{z+1}{kz-1}
\]
maps the domain $R_2$ onto the unbounded doubly connected domain $R_3$ exterior to the two slits $[-1,0]$ and $[s,+\infty]$ where $s=\Psi(1)=4k/(1-k)^2$.
Thus, the capacity of the domain $R_2$ equals to the capacity of $R_3$ which can be expressed by~\cite[5.60 (1)]{vu88}
\begin{equation}\label{eq:t-b-cap}
\capa(R_2)=\frac{\pi}{\mu(1/\sqrt{1+s})}, \quad s=\frac{4k}{(1-k)^2}.
\end{equation}
Here $\mu$ is the function defined in \eqref{mu}.

By~\cite[119.03]{bf}, the domain $R_2$ can be mapped conformally also onto the unbounded doubly connected domain $R_4$ exterior to the two slits $[-t/2-\i b/2,-t/2+\i b/2]$ and $[t/2-\i b/2,t/2+\i b/2]$ with
\[
t=\frac{2}{k}\left(E(k)-(1-k^2a^2)K(k)\right), \quad
b=\frac{2}{k}\left(E(k'_1,k')-k^2a^2F(k'_1,k'))\right)
\]
where the functions $E(k), K(k)$ are defined in \eqref{defE} and \eqref{defK}, resp., and
\[
F(z,k)=\int_0^z \frac{dw}{\sqrt{(1-w^2)(1-k^2w^2)}}, \quad
E(z,k)=\int_0^z \sqrt{\frac{1-k^2w^2}{1-w^2}}dw,
\]
and
\[
k'=\sqrt{1-k^2}, \quad a=\frac{E(k')}{k^2K(k')} \quad
k_1=\frac{k}{k'}\sqrt{a^2-1}, \quad k'_1=\sqrt{1-k_1^2}.
\]
Hence $\capa(R_4)=\capa(R_2)$.
Further, it is clear that the domain $R_1$ can be conformally mapped by the function $\hat\Psi(z)=\i(z-1/2)$ onto the domain $R_4$ if we choose $k$ such that $t=b=1$. Thus, the exact capacity of $R_1$ is given by
\begin{equation}\label{eq:cap-R1-v}
\capa(R_1)=\frac{\pi}{\mu(1/\sqrt{1+s})}, \quad s=\frac{4k}{(1-k)^2}
\end{equation}
where $k$ satisfies the equations
\begin{equation}\label{eq:t-b}
1=\frac{2}{k}\left(E(k)-(1-k^2a^2)K(k)\right), \quad
1=\frac{2}{k}\left(E(k'_1,k')-k^2a^2F(k'_1,k'))\right).
\end{equation}

The equations~\eqref{eq:t-b} are solved using Mathematica for $k$ and the value of the capacity of $R_1$ computed through~\eqref{eq:cap-R1-v} is $2.1157789709245134$. This value agrees with the value presented at the bottom of Table~\ref{tab:two-rec} with relative error $1.5\times10^{-9}$.

\subsection{A vertical rectangle in the upper half-plane}

\begin{figure}[ht] %
\centerline{
\hfill
\scalebox{0.35}{\includegraphics[trim=0 0 0 0cm,clip]{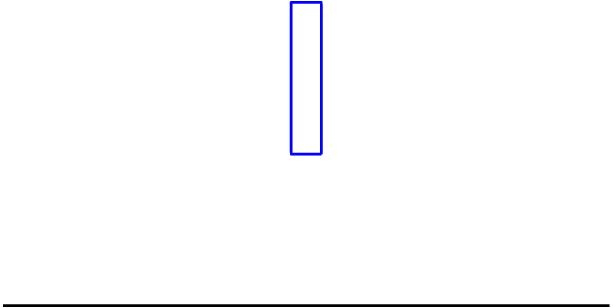}}
\hfill
\scalebox{0.35}{\includegraphics[trim=0 0 0 0cm,clip]{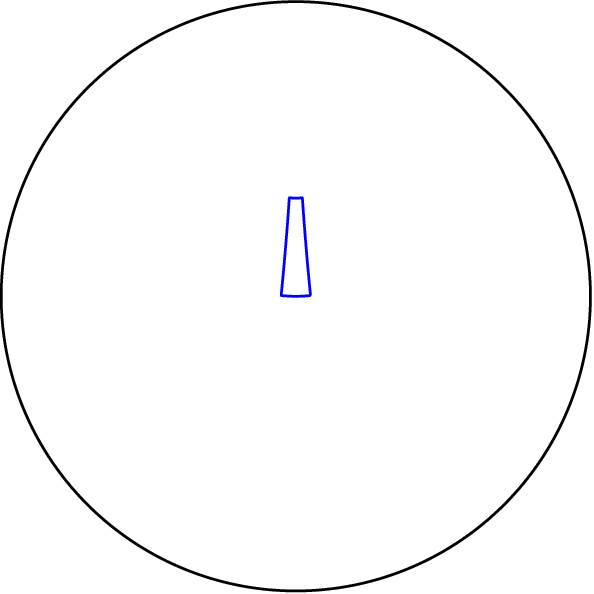}}
\hfill
}
\caption{The domain $G$ in the exterior of a vertical rectangular closed set in the upper half-plane (left) and its image $\hat G$ under the auxiliary map $\Psi$ (right) for $d=0.1$.}
\label{fig:rec-half}
\end{figure}

Consider the doubly connected domain $G$ exterior to the rectangular closed set
 \[
 [0.5-d,0.5+d] \times [1,2]
 \]
in the upper half-plane where $0<d<0.5$ (see the Figure~\ref{fig:rec-half} for $d=0.1$).
The auxiliary map
\begin{equation}\label{eq:half-rec-psi}
w=\Psi(z)=\frac{\i z+1}{z+\i}
\end{equation}
is used to transform the domain $G$ onto a domain $\hat G$ interior to the unit disk and exterior to the piecewise smooth Jordan curve $L$ which is the image of the rectangle under the map $\Psi$. Then $G$ and $\hat G$ have the same capacities. We use the function \verb|annq| with $n=2^{15}$ to compute the capacity of $\hat G$ for several values of $d$. When $d=0$, the rectangle reduced to the slit $[\i,2\i]$. For the upper half-plane with the slit $[\i,2\i]$, we can use the numerical method presented in Section~\ref{sec:half-plane} to compute the capacity of the domain exterior to this slit in the upper half-plane. The results are presented in Table~\ref{tab:rec-half}. The exact value of the capacity of the domain exterior to slit  $[\i,2\i]$ in the upper half-plane can be computed from~\eqref{eq:cap-half-seg} and is equal to $2\pi/\mu(1/3)=2.55852314234201$. The result presented at the bottom of Table~\ref{tab:rec-half} agrees with the exact value with relative error $4.7\times10^{-13}$.

\begin{table}[ht]
\caption{The capacity values for the domain exterior to a vertical rectangle in the upper half-plane.}
\label{tab:rec-half}%
\begin{tabular}{l|l|l}\hline
 $d$  & $\capa(G)$        & Time (sec)\\ \hline
0.4   & 3.71752232703208  & 2.4  \\
0.3		& 3.46693660197964	& 2.2  \\
0.2   & 3.20488821317939  & 2.2  \\
0.1   & 2.9209225535743   & 2.3  \\
0.05  & 2.76128813737089  & 2.6  \\
0.02  & 2.65173985860514  & 3.0  \\
0.01  & 2.60986001541974  & 3.9  \\
0.005 & 2.58658944233183  & 5.4  \\
0     & 2.55852314234082  & 16.4  \\
\hline
\end{tabular}
\end{table}

\subsection{A horizontal rectangle in the upper half-plane}

\begin{figure}[ht] %
\centerline{
\hfill
\scalebox{0.35}{\includegraphics[trim=0 0 0 0cm,clip]{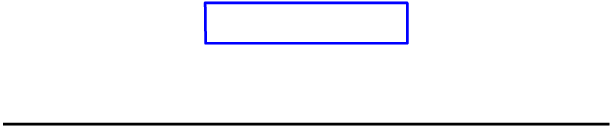}}
\hfill
\scalebox{0.35}{\includegraphics[trim=0 0 0 0cm,clip]{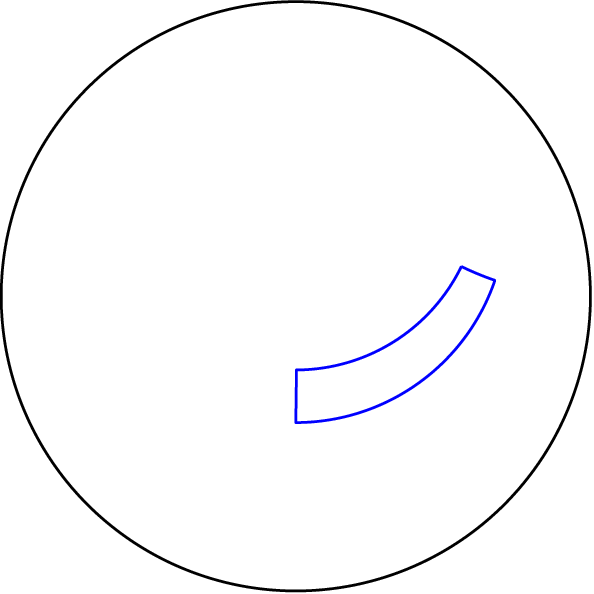}}
\hfill
}
\caption{The domain $G$ in the exterior of a horizontal rectangular closed set in the upper half-plane (left) and its image $\hat G$ under the auxiliary map $\Psi$ (right) for $d=0.1$.}
\label{fig:rec-half-h}
\end{figure}

Consider the doubly connected domain $G$ exterior to the rectangular closed set
 \[
 [0,1]\times[0.5-d,0.5+d]
 \]
in the upper half-plane where $0<d<0.5$ (see the Figure~\ref{fig:rec-half-h}). By symmetry, the capacity for this domain is $2$ times the capacity for the two rectangles case considered in Subsection~\ref{sec:two-rec}.

As in the previous example, the the auxiliary map $w=\Psi(z)$ in~\eqref{eq:half-rec-psi} is used to transform the domain $G$ onto a domain $\hat G$ interior to the unit disk and exterior to a piecewise smooth Jordan curve (see Figure~\ref{fig:rec-half-h}). Then $G$ and $\hat G$ have the same capacities. We use the function \verb|annq| with $n=2^{15}$ to compute the capacity of $\hat G$ for several values of $d$.
When $d=0$, the rectangle reduced to the slit $[0.5\i,1+0.5\i]$. By symmetry, the capacity for the half-plane with the horizontal slit $[0.5\i,1+0.5\i]$ is $2$ times the capacity for the cases of the domain exterior to the two horizontal slits $[\i/2,1+\i/2]$ and $[-\i/2,1-\i/2]$ considered in Subsection~\ref{sec:two-rec}. Thus, according to the exact capacity presented in Subsection~\ref{sec:two-rec}, the exact capacity for the upper half-plane with the horizontal slit $[0.5\i,1+0.5\i]$ is $4.23155794184903$.

For numerical computing of the capacity of the upper half-plane with the slit $[0.5\i,1+0.\i]$, we use the method described in Section~\ref{sec:half-plane}. The obtained result is presented at the bottom of Table~\ref{tab:rec-half-h}. The computed approximate value agrees with the exact value with relative error $1.1\times10^{-9}$.

Finally, the third column in Table~\ref{tab:rec-half-h} shows halves of the computed values of the capacity for the domain presented in this section. The values presented in the third column agrees with the results presented in Table~\ref{tab:two-rec} for two rectangle case.

\begin{table}[ht]
\caption{The capacity values for the domain exterior to a horizontal rectangle in the upper half-plane.}
\label{tab:rec-half-h}%
\begin{tabular}{l|l|l|l}\hline
 $d$  & $\capa(G)$        & $\capa(G)/2$        & Time (sec)\\ \hline
0.4   & 15.1134561077006  & 7.5567280538503     & 2.6  \\
0.3		& 9.10569023215289	& 4.55284511607644    & 2.5  \\
0.2   & 6.77138475734822  & 3.38569237867411    & 2.3  \\
0.1   & 5.37377572427995  & 2.68688786213998    & 2.4  \\
0.05  & 4.81109439601605  & 2.40554719800803    & 2.5  \\
0.02  & 4.48126118775531  & 2.24063059387766    & 3.1  \\
0.01  & 4.36525097360269  & 2.18262548680134    & 4.1  \\
0.005 & 4.30323272661648  & 2.15161636330824    & 5.7  \\
0     & 4.2315579463472   & 2.1157789731736     & 23.0  \\
\hline
\end{tabular}
\end{table}


\section{The hyperbolic capacity and the elliptic capacity}

Let $E$ be a compact and connected set (not a single point) in the unit disk $\D$.
In this section, we use the MATLAB function \verb|annq| in Subsection~\ref{sc:comp-cap} to compute the hyperbolic capacity and the elliptic capacity of the set $E$. Both the hyperbolic capacity and the elliptic capacity are invariants under conformal mappings.

\subsection{The hyperbolic capacity}

The hyperbolic capacity of $E$, $\caph(E)$, is defined by~\cite[p.~19]{Vas02}
\begin{equation}\label{eq:hyp-cap}
\caph(E)=\lim_{n\to\infty}\left[\max_{z_1,\ldots,z_n\in E}\prod_{1\le k<j\le n}
\left|\frac{z_k-z_j}{1-z_k\overline{z_j}}\right|
\right]^{\frac{2}{n(n-1)}}.
\end{equation}
For the hyperbolic capacity, we assume $G$ is the bounded doubly connected domain defined by $G=\D\backslash E$ such that $G$ can be mapped conformally onto an annulus $q<|w|<1$. Hence the hyperbolic capacity $\caph(E)$ is given by~\cite{DuKu}
\begin{equation}\label{eq:hyp-cap3}
\caph(E)=q.
\end{equation}
The constant $q$ can be computed by the function \verb|annq|.

\subsection{The elliptic capacity}

For the compact and connected set $E$, we define the antipodal set $E^\ast=\{-1/\overline{a}\,:\,a\in E\}$. Since we assume $E\subset\D$, we have $E\cap E^\ast=\emptyset$ (in this case, the set $E$ is called ``elliptically schlicht''~\cite{DuKu}). The elliptic capacity of $E$, $\cape(E)$, is defined by~\cite{DuKu}
\begin{equation}\label{eq:ell-cap}
\cape(E)=\lim_{n\to\infty}\left[\max_{z_1,\ldots,z_n\in E}\prod_{1\le k<j\le n}
\left|\frac{z_k-z_j}{1+z_k\overline{z_j}}\right|
\right]^{\frac{2}{n(n-1)}}.
\end{equation}

To compute the elliptic capacity, we assume $G$ is the doubly connected domain between $E$ and $E^\ast$ such that $G$ can be mapped conformally onto an annulus $r<|w|<1/r$. Then the elliptic capacity is given by~\cite{DuKu}
\[
\cape(E)=r.
\]
Here, the domain $G$ could be bounded or unbounded. We shall use the method described in Section~\ref{sec-cm} to map the domain $G$ onto an annulus $q<|w|<1$ which is conformally equivalent to the annulus $r<|w|<1/r$ with $r=\sqrt{q}$. Thus, we have
\begin{equation}\label{eq:ell-cap2}
\cape(E)=\sqrt{q}.
\end{equation}
We compute $q$ using the function \verb|annq|.

Finally, as our interest in this paper is only in closed and connected subsets $E$ of the unit disk $\D$ and comparing numerically between the values of $\cape(E)$ and $\caph(E)$, it is worth mentioning that Duren and K\"uhnau~\cite{DuKu} have proved that
\[
\cape(E)\le\caph(E),
\]
with equality if and only if $E=-E$. This inequality is verified numerically in the following numerical examples.

\subsection{A disk}

As our first example, we compute the hyperbolic capacity and the elliptic capacity of the disk $E=\{z\in\C\;:\;|z|\le r\}$, $0<r<1$. For this set $E$, both capacities are equal where~\cite{kir,DuKu}
\[
\caph(E)=\cape(E)=r.
\]

\begin{figure}[ht] %
\centerline{
\hfill
\scalebox{0.35}{\includegraphics[trim=0cm 0cm 0cm 0cm,clip]{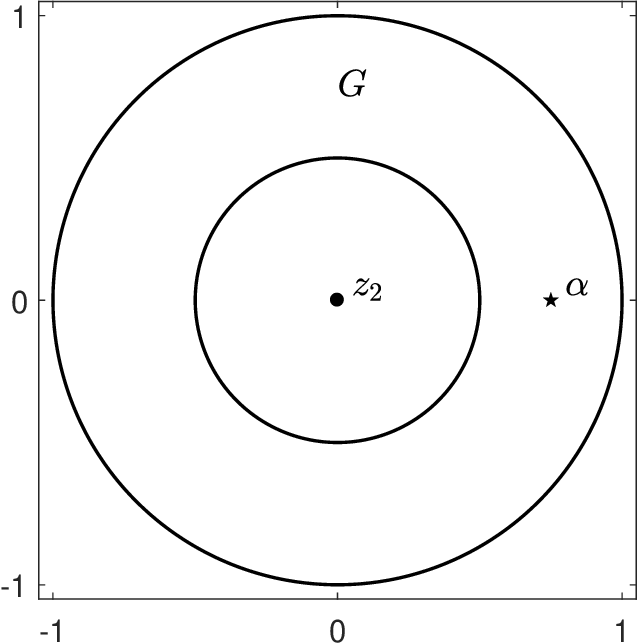}}
\hfill
\scalebox{0.35}{\includegraphics[trim=0cm 0cm 0cm 0cm,clip]{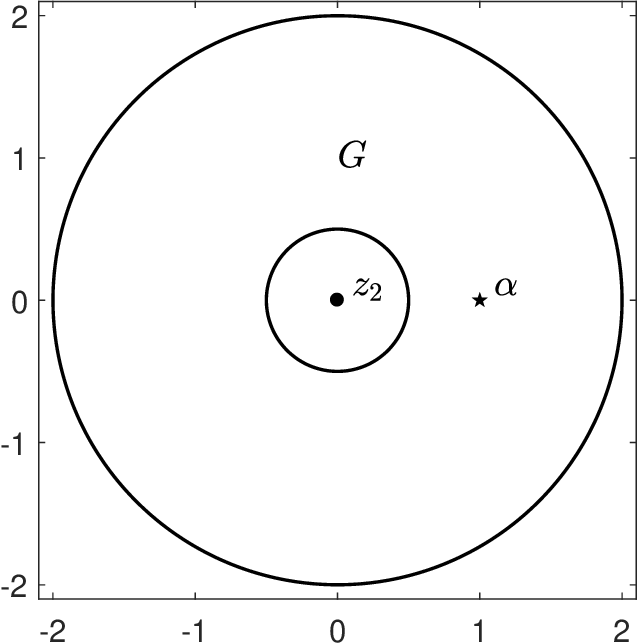}}
\hfill
}
\caption{The domain $G$ for computing the hyperbolic capacity (left) and the elliptic capacity (right) of $E=\{z\in\C\;:\;|z|\le r\}$ for $r=0.5$.}
\label{fig:LHEcap_circ}
\end{figure}

\begin{figure}[ht] %
\centerline{\scalebox{0.4}{\includegraphics[trim=0cm 0cm 0cm 0cm,clip]{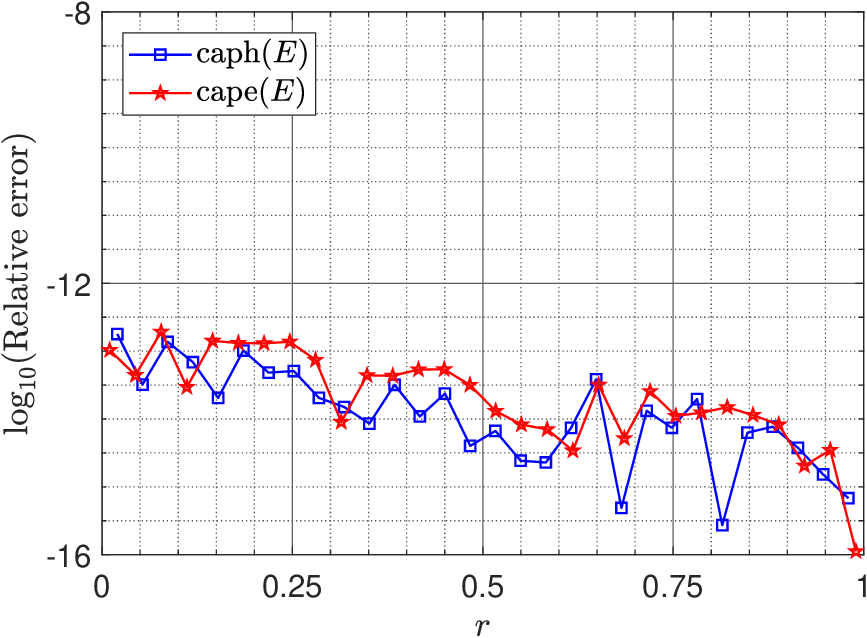}}}
\caption{The relative error in the computed capacities $\cc(E)$, $\caph(E)$, and $\cape(E)$ for the disk  $E=\{z\in\C\;:\;|z|\le r\}$.}
\label{fig:LHEcap_cir}
\end{figure}

For computing $\caph(E)$, we use the function \verb|annq| with $\alpha=(1+r)/2$ and $z_2=0$ to compute the value of $q$ for the conformal map of the doubly connected domain $G=\D\backslash E$ (see Figure~\ref{fig:LHEcap_circ} (left)) onto the annulus $q<|w|<1$ and hence $\caph(E)=q$. For $\cape(E)$, the domain $G$ between $E$ and $E^\ast$ is the bounded doubly connected domain $r<|z|<1/r$ (see Figure~\ref{fig:LHEcap_circ} (right)). We use the MATLAB function \verb|annq| with $\alpha=1$ and $z_2=0$ to compute the value of $q$ for the conformal map of this domain $G$ onto the annulus $q<|w|<1$ and hence $\cape(E)=\sqrt{q}$. For both cases, we use $n=2^{12}$ and $0.02\le r\le 0.98$.
The relative error in the obtained results for $\caph(E)$ and $\cape(E)$ are shown in Figure~\ref{fig:LHEcap_cir}.

\subsection{A square}

\begin{figure}[ht] %
\centerline{
\hfill
\scalebox{0.35}{\includegraphics[trim=0cm 0cm 0cm 0cm,clip]{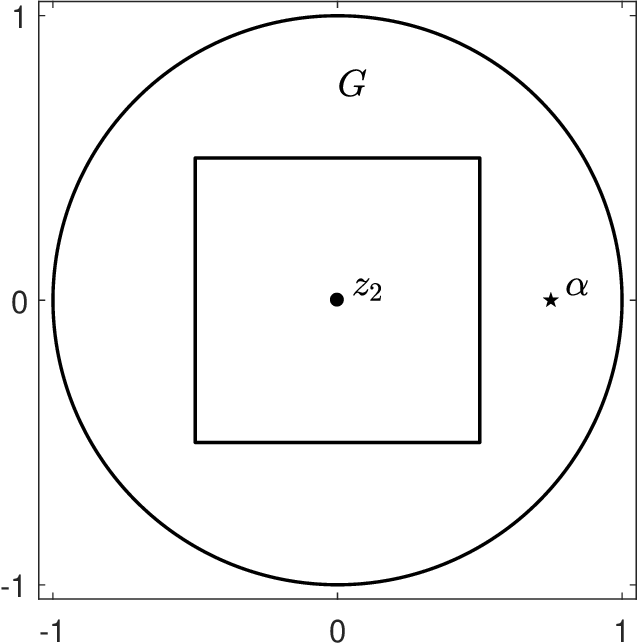}}
\hfill
\scalebox{0.35}{\includegraphics[trim=0cm 0cm 0cm 0cm,clip]{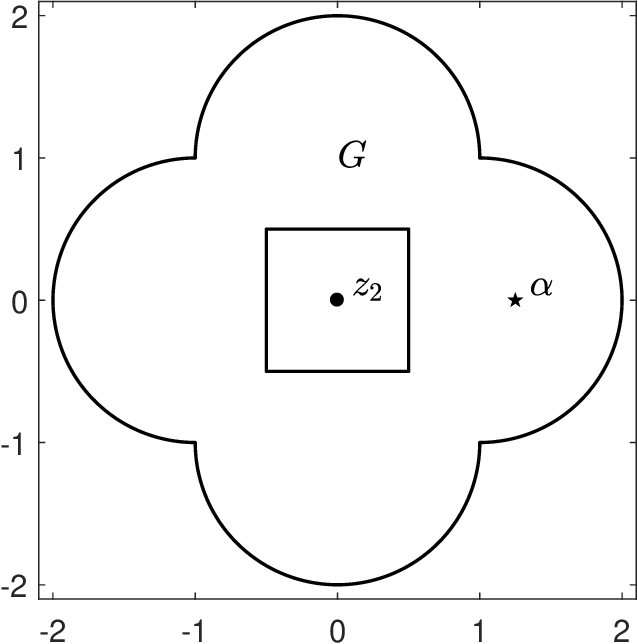}}
\hfill
}
\caption{The domain $G$ for computing the hyperbolic capacity (left) and the elliptic capacity (right) of $E=[-r,r]\times[-r,r]$ for $r=0.5$.}
\label{fig:LHEcap_sqc}
\end{figure}

\begin{figure}[!h] %
\centerline{\scalebox{0.4}{\includegraphics[trim=0cm 0cm 0cm 0cm,clip]{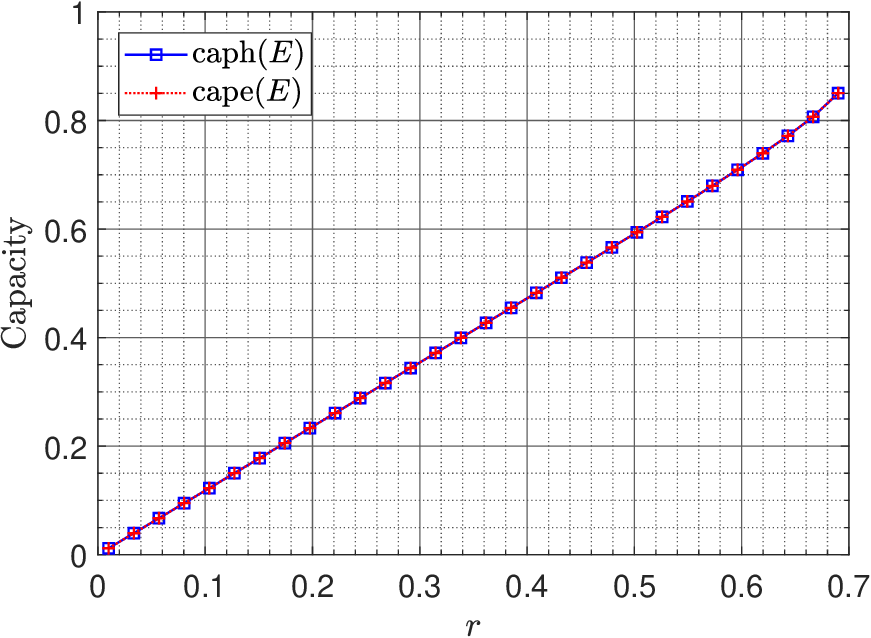}}}
\caption{The capacities $\caph(E)$ and $\cape(E)$ for $E=[-r,r]\times[-r,r]$.}
\label{fig:LHEcap_sq}
\end{figure}

For the second example, we assume $E$ is the closed set $[-r,r]\times[-r,r]$, $0<r<1/\sqrt{2}$.
For computing $\caph(E)$, the domain $G$ is the bounded doubly connected domain exterior to the square and interior to the unit circle (see Figure~\ref{fig:LHEcap_sqc} (left)). We use the function \verb|annq| with $\alpha=(1+r)/2$ and $z_2=0$ to compute $q$ and then $\caph(E)=q$. 
For $\cape(E)$, the domain $G$ is the bounded doubly connected domain between $E$ and $E^\ast$  (see Figure~\ref{fig:LHEcap_sqc} (right)). Hence, $\cape(E)=\sqrt{q}$ where $q$ is computed using the function \verb|annq| with $\alpha=(r+1/r)/2$ and $z_2=0$.  
For both cases, we use $n=2^{13}$ for $0.02\le r\le 0.69$. The obtained results are shown in Figure~\ref{fig:LHEcap_sq}. This set $E$ is symmetric where $E=-E$, and hence $\caph(E)=\cape(E)$. 

\subsection{Amoeba-shaped boundary}

\begin{figure}[ht] %
\centerline{
\hfill
\scalebox{0.35}{\includegraphics[trim=0cm 0cm 0cm 0cm,clip]{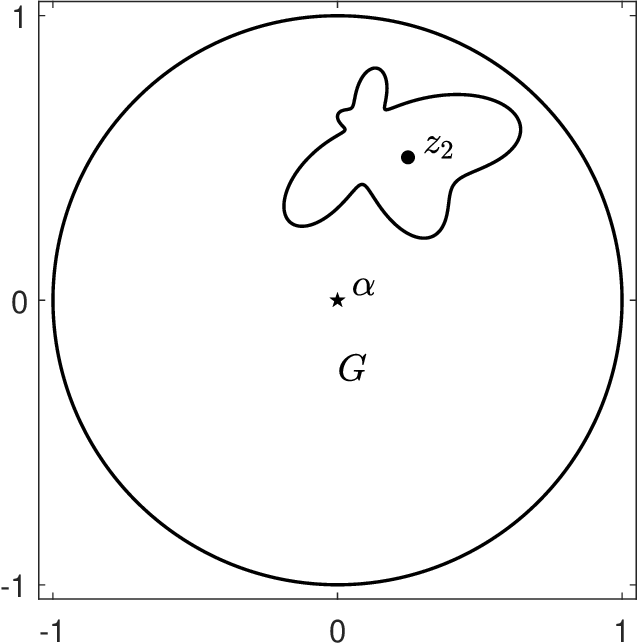}}
\hfill
\scalebox{0.35}{\includegraphics[trim=0cm 0cm 0cm 0cm,clip]{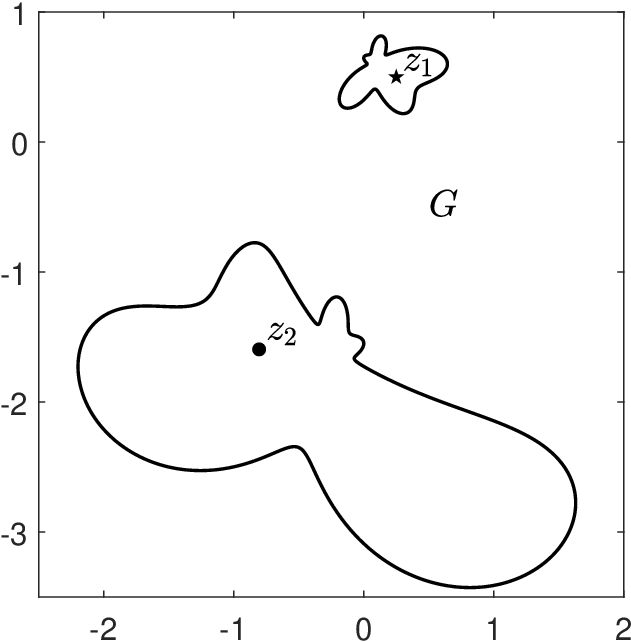}}
\hfill
}
\caption{The domain $G$ for computing the hyperbolic capacity (left) and the elliptic capacity (right) of the closed region bordered by the amoeba-shaped boundary $E$.}
\label{fig:LHEcap_amc}
\end{figure}

For the third example, we compute $\caph(E)$ and $\cape(E)$ of $E$ where $E$ is the closed region bordered by the amoeba-shaped boundary $L$ with the parametrization
\[
\eta(t)=0.1+0.6\i+0.2\left(e^{\cos t}\cos^22t+e^{\sin t}\sin^22t\right)e^{-\i t}, \quad 0\le t\le 2\pi.
\]
For the hyperbolic capacity $\caph(E)$, the domain $G$ is the bounded doubly connected domain exterior to the curve $L$ and interior to the unit circle (see Figure~\ref{fig:LHEcap_amc} (left)). Then $\caph(E)=q$ where $q$ is computed using the function \verb|annq| with $\alpha=0$ and $z_2=0.25+0.5\i$. To compute $\cape(E)$, the domain $G$ is the unbounded doubly connected domain exterior to $E$ and $E^\ast$ (see Figure~\ref{fig:LHEcap_amc} (right)). We use the function \verb|annq| with $z_1=0.25+0.5\i$ and $z_2=-1/\overline{z_1}$ to compute the value of $q$ and hence $\cape(E)=\sqrt{q}$.
The approximate values of the capacities $\caph(E)$ and $\cape(E)$ for several values of $n$ are shown in Table~\ref{tab:am}. As the set $E$ is not symmetric, the presented numerical results confirmed the inequality $\cape(E)<\caph(E)$.

\begin{table}[!h]
\caption{The approximate values of the capacities of the closed region bordered by the amoeba-shaped boundary $E$.}
\label{tab:am}%
\begin{tabular}{l|l|l}\hline
   $n$  & $\caph(E)$         & $\cape(E)$        \\ \hline
    64  &0.521349946390291   & 0.25872431985379  \\
   128  &0.521358819409768   &0.258724285703159  \\
   256  &0.521358832558364   &0.258724285703154  \\
   512  &0.521358832558375   &0.258724285703153  \\
  1024  &0.52135883255838    &0.258724285703155  \\
  2048  &0.521358832558369   &0.258724285703154  \\
  4096  &0.521358832558378   &0.258724285703156  \\
\hline %
\end{tabular}
\end{table}


\section{Concluding Remarks}

Conformal invariants are important tools for complex analysis with many applications. However, these invariants can be expressed explicitly only in very few special cases. Thus, numerical methods are required to compute these invariants. 
A numerical method for computing some conformal invariants is presented in this paper. The method can be used for domains with different types of boundaries including domains with smooth or piecewise smooth boundaries. The performance and the accuracy of the presented method is compared to analytic solutions or to previous results whenever analytic solutions or previous results are available. 
Further, a MATLAB implementation of the proposed method is given in the MATLAB function \verb|annq| in Subsection~\ref{sc:comp-cap}. This MATLAB function was used in almost all examples in this paper to compute the conformal capacity, the hyperbolic capacity and the elliptic capacity. For some examples, an auxiliary procedure is required before using the function \verb|annq|. The computer codes of the presented computations can be found in the link \url{https://github.com/mmsnasser/cci}.

\bibliographystyle{siam}


\end{document}